\documentclass[12pt]{article}

\usepackage{amsfonts, amsmath, amssymb, amsthm, mathtools}

\usepackage[margin=1.1in]{geometry}

\usepackage{graphicx}
\usepackage{float}
\usepackage{booktabs}
\usepackage{longtable}
\usepackage{subcaption}
\usepackage{array}
\usepackage{multirow}
\usepackage{threeparttable}

\usepackage{blkarray}

\usepackage[numbers]{natbib}
\usepackage{url}

\usepackage{enumitem}
\usepackage{setspace}
\usepackage{microtype}

\usepackage{pdflscape}

\usepackage{tikz}

\newcommand\hlight[1]{%
  \tikz[overlay,remember picture,
        baseline=-\the\dimexpr\fontdimen22\textfont2\relax]
  \node[rectangle,
        fill=blue!50,
        rounded corners,
        fill opacity=0.2,
        draw,
        thick,
        text opacity=1] {$#1$};%
}

\newtheorem{theorem}{Theorem}[section]

\newtheorem{lemma}[theorem]{Lemma}

\theoremstyle{definition}
\newtheorem{example}[theorem]{Example}

\theoremstyle{remark}
\newtheorem{remark}[theorem]{Remark}

\DeclareMathOperator{\tr}{tr}
\DeclareMathOperator{\diag}{diag}

\begin{document}

\singlespace

\title{Counting five-node subgraphs}

\author{\textsc{Steve Lawford}\footnote{ENAC (University of Toulouse), 7 avenue Edouard Belin, BP 54005, 31055, Toulouse, Cedex 4, France (email: steve.lawford@enac.fr).}
}

\date{}

\singlespace

\maketitle

\vspace{0.5cm}

\begin{abstract}

\noindent For a simple undirected graph $G$, we derive exact formulae for the number of copies of each of the 21 non-isomorphic connected graphs on five vertices. For nine of these graphs (the stingray, spinning top, kite, ufo, crown, envelope, lamp, arrowhead, and cat's cradle) the resulting analytical formulae appear to provide substantially new formulations; we also identify and correct five erroneous published formulae for the 5-path, banner, and lollipop. The proofs use elementary combinatorial arguments, organized around recurring constructions based on incident structures, walks, common neighbourhoods, and neighbourhood subgraphs, several of which extend naturally to larger subgraphs. We illustrate the formulae by deriving analytical results for several regular graphs and by applying them to a real-world network, where induced five-node subgraph counts, obtained as linear combinations of general counts, are compared with Erd\H{o}s--R\'enyi and degree-preserving null ensembles.

\end{abstract}

\noindent\textbf{Keywords:} subgraph counting; graph enumeration; complex networks; small subgraphs.

\medskip

\newpage

\section{Introduction}

Counts of small connected subgraphs are fundamental descriptors of
network structure. They are used, for example, for network classification and to assess the statistical significance of small topological structures that occur frequently in real-world networks, where they may have specific functional roles \cite{alon07, benson_etal16, estrada07, hayes_etal13, itzkovitz_alon05, kashtan_etal04, kuchaiev_etal10, milenkovic_przulj08, milo_etal02,
xia_etal19}. Subgraph-count-based statistics are also used for network comparison and model assessment \cite{ospina-forero_etal19}, while systematic small-subgraph analysis has been extended to multilayer networks \cite{sallmen_etal22}. Their enumeration is closely related to the classical subgraph isomorphism problem in theoretical computer science. Exact analytical count formulae have been studied since at least the early 1970s \cite{harary_manvel71}, but their application becomes increasingly costly as the number of nodes in the target subgraph grows. This combinatorial explosion has motivated extensive work on exact enumeration and approximate sampling algorithms, with applications to increasingly large datasets (an early subgraph-counting algorithm is given in \cite{itai_rodeh78}; for recent work, see references in \cite{bera_etal20} and the comprehensive survey by
\cite{ribeiro_etal21}). In contrast, there are almost no complete
analytical treatments of exact subgraph counting on five nodes, despite the theoretical insights that can be gained from such formulae and their usefulness in applications such as
\cite{agasse-duval_lawford25, allen-perkins_etal17, lawford_mehmeti20}.\\
\\
\noindent We fill this gap by providing exact count formulae for all 21 non-isomorphic connected subgraphs on five nodes (Table
\ref{tab:subgraphs5}) in simple unweighted undirected graphs. The
formulae are derived using elementary combinatorial arguments organized around four recurring constructions. For nine subgraphs, the resulting analytical formulae appear to provide substantially new formulations; we also identify and correct five erroneous published formulae for the 5-path, banner, and lollipop. The neighbourhood construction can be formalized as a general lifting identity that gives a common derivation of several formulae in Theorem~\ref{thm:analytic_count_5node} and extends directly to suitable larger subgraphs. Exact results for the eight connected subgraphs on three or four nodes are well known
\cite{agasse-duval_lawford25, estrada11, estrada_knight15}. The number of non-isomorphic connected graphs nevertheless grows rapidly:
there are already 112 on six nodes, making a complete analytical treatment at that order impractical with the present methods, although individual larger subgraphs can still be treated. See A001349 of the Online Encyclopedia of Integer Sequences (\url{https://oeis.org/A001349}) for the number of non-isomorphic connected graphs on $n$ nodes.\\
\\
\noindent The paper is organized as follows. In Section
\ref{sec:main_result} we present the main result and the four recurring strategies used in its proof. In Section \ref{sec:examples} we specialize selected count formulae to three classical regular graphs: the complete graph $K_n$, the complete $N$-partite graph
$K_{n_1,n_2,\ldots,n_N}$, and the regular ring lattice. These examples
develop combinatorial intuition for the formulae. We then show how
alternative derivations can lead to quite different formulae and,
sometimes, to error, and illustrate the implementation of all 21 count
formulae on a real-world network, including induced-count comparisons
with Erd\H{o}s--R\'enyi and degree-preserving null ensembles. Section
\ref{sec:conclusions} concludes. Appendix \ref{sec:induced_counts}
discusses the well-known result that induced subgraph counts can be
obtained as linear combinations of general subgraph counts.

\clearpage

\begin{landscape}

\vspace*{\fill}

\begin{center}

\resizebox{0.92\linewidth}{!}{%
\begin{tabular}{ccccccc}
\includegraphics[height=1in]{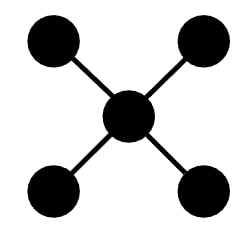} &
\includegraphics[height=1in]{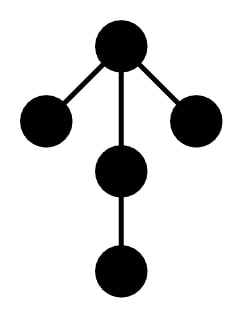} &
\includegraphics[height=1in]{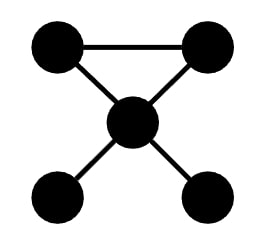} &
\includegraphics[height=1in]{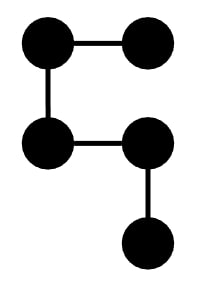} &
\includegraphics[height=1in]{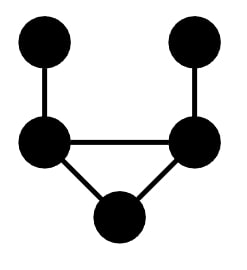} &
\includegraphics[height=1in]{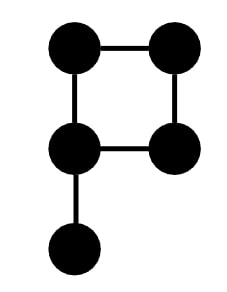} &
\includegraphics[height=1in]{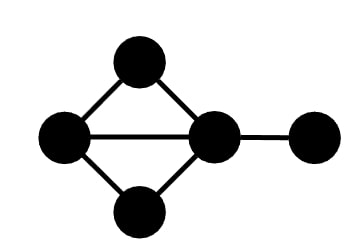} \\

5-star $M_{75}^{(5)}$ &
5-arrow $M_{77}^{(5)}$ &
cricket $M_{79}^{(5)}$ &
5-path $M_{86}^{(5)}$ &
bull $M_{87}^{(5)}$ &
banner $M_{94}^{(5)}$ &
stingray $M_{95}^{(5)}$\\[1ex]

\includegraphics[height=1in]{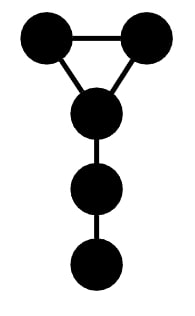} &
\includegraphics[height=1in]{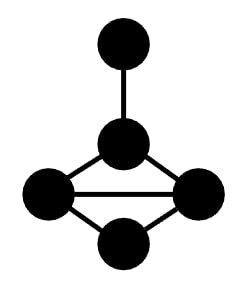} &
\includegraphics[height=1in]{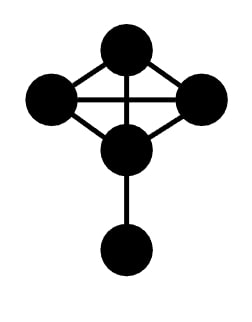} &
\includegraphics[height=1in]{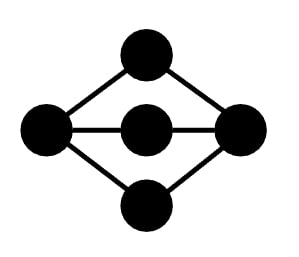} &
\includegraphics[height=1in]{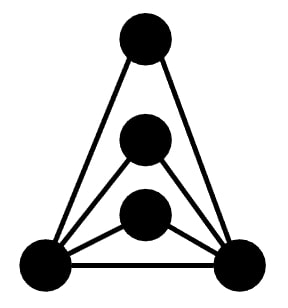} &
\includegraphics[height=1in]{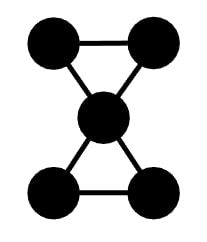} &
\includegraphics[height=1in]{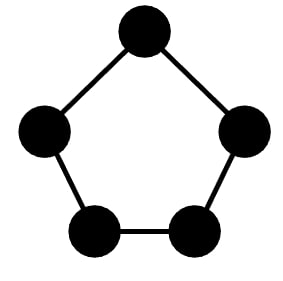} \\

lollipop $M_{117}^{(5)}$ &
spinning top $M_{119}^{(5)}$ &
kite $M_{127}^{(5)}$ &
ufo $M_{222}^{(5)}$ &
chevron $M_{223}^{(5)}$ &
hourglass $M_{235}^{(5)}$ &
5-circle $M_{236}^{(5)}$\\[1ex]

\includegraphics[height=1in]{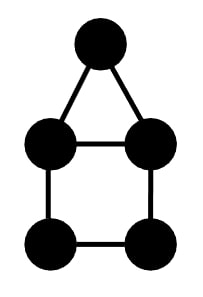} &
\includegraphics[width=1in]{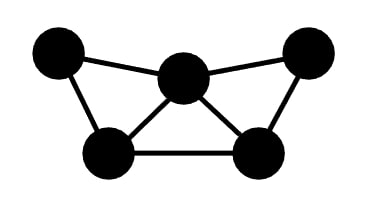} &
\includegraphics[height=1in]{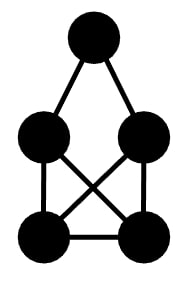} &
\includegraphics[height=1in]{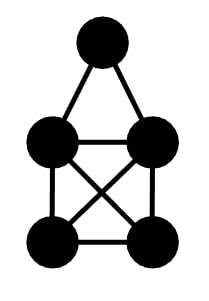} &
\includegraphics[height=1in]{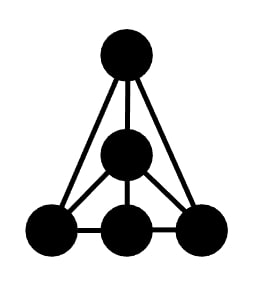} &
\includegraphics[height=1in]{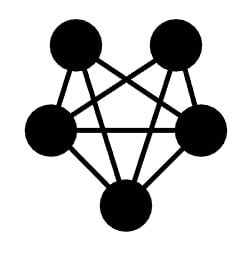} &
\includegraphics[height=1in]{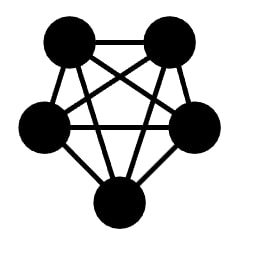} \\

house $M_{237}^{(5)}$ &
crown $M_{239}^{(5)}$ &
envelope $M_{254}^{(5)}$ &
lamp $M_{255}^{(5)}$ &
arrowhead $M_{507}^{(5)}$ &
cat's cradle $M_{511}^{(5)}$ &
5-complete $M_{1023}^{(5)}$
\end{tabular}%
}

\captionof{table}{The 21 non-isomorphic connected subgraphs on five
nodes, denoted $M_{a}^{(b)}$ (see Section~\ref{sec:notation} for notation).}
\label{tab:subgraphs5}

\end{center}

\vspace*{\fill}

\end{landscape}

\subsection{Historical context}\label{sec:history}
Subgraph enumeration has a long history in combinatorics, with early
work by Harary and Manvel in 1971 \cite{harary_manvel71} deriving
exact formulae for cycles of small length; one of the earliest papers
in the field appears to be by Cartwright and Gleason in 1966
\cite{cartwright_gleason66}. A major advance followed in 1997 with the work of Noga Alon, Yuster, and Zwick \cite{alon_etal97}, who developed efficient methods for finding and counting fixed-length cycles. Since then, progress on analytical enumeration has been limited, with most advances focused on algorithmic methods or induced subgraphs (and motifs). In 2002, Milo, Shen-Orr, Itzkovitz, Kashtan, Chklovskii and Uri Alon \cite{milo_etal02} introduced the concept of network motifs as recurring small induced subgraphs in real networks. This sparked renewed interest in systematically counting all small subgraphs of a certain size. A key insight was that induced subgraph counts can be derived as linear combinations of general subgraph counts, via inclusion–exclusion on the presence of extra edges. This principle essentially relies on the M{\"o}bius inversion on the lattice of subgraphs. From 2011 onwards, analytic subgraph counting started to appear in important textbook treatments of graph theory and network science e.g. \cite{estrada11, estrada_knight15}.\\
\\
Related analytical work has exploited additional structure in regular
graphs. Beezer \cite{beezer01} studies linear relations among counts of isomorphism classes of subgraphs in regular graphs, within a more general framework for configurations in combinatorial designs. Minchenko and Wanless \cite{minchenko_wanless14} use generating functions for closed walks to express spectral moments of regular graphs as linear combinations of counts of certain connected subgraphs. These approaches are complementary to the present work: they establish structural relations among selected subgraph counts in regular host graphs, whereas our aim is to derive explicit count formulae for all connected five-node subgraphs in arbitrary simple undirected graphs. A recent survey of motif-discovery methods discusses exact and approximate frequency algorithms, including methods for graphlet mining \cite{jazayeri_yang20}. Related work continues to develop scalable methods for subgraph enumeration, including distributed enumeration of small directed subgraphs \cite{levinas_etal22} and memory-efficient methods based on local subgraph statistics in very large networks \cite{weedage_etal21}. Recent developments in the theoretical computer science community have focused on algorithmic frameworks (e.g., ESCAPE \cite{pinar_etal17}, EVOKE \cite{pashanasangi_seshadhri20}, HUGE \cite{yang_etal21}) that can count all five-node subgraphs, but few works have produced closed-form formulae; see also \cite{conte_etal22, klusowski_wu18, kowaluk_etal13}. Bera et al. \cite{bera_etal20} describe a ``chasm at size six'', where
techniques effective up to five nodes become infeasible without new ideas.

\subsection{Notation and preliminaries}\label{sec:notation}
Let $G=(V,E)$ be a graph with node set $V$ and edge set $E$. We write
$n=|V|$ and $m=|E|$ for the numbers of nodes and edges of $G$. Let $g$
be the $n\times n$ adjacency matrix of $G$, with representative element $g_{ij}$, which takes value one when an edge is present between nodes $i$ and $j$, and zero otherwise. We use $\{i,j\}\in E$ to denote an undirected edge between nodes $i$ and $j$, and say that the nodes are directly connected. A graph is simple and unweighted if $g_{ii}=0$ and $g_{ij}\in\{0,1\}$, and undirected if $g_{ij}=g_{ji}$. When an orientation of an undirected edge is required, we write $(i,j)$ or $(j,i)$ for its two orientations. A walk from node $i$ to node $j$ is a sequence of nodes $i_1,\ldots,i_{R+1}$ with $i_1=i$,
$i_{R+1}=j$, and $\{i_r,i_{r+1}\}\in E$ for $r=1,\ldots,R$. A path is a walk whose nodes are distinct. A graph is connected if there is at least one path between every pair of nodes $i$ and $j$.\\
\\
\noindent
We use $\Gamma_G(i)=\{j:\{i,j\}\in E\}$ to denote the neighbourhood of node $i$ in $G$, with cardinality equal to its degree $k_i=|\Gamma_G(i)|=\sum_j g_{ij}$. Let $P(k)$ denote the distribution of node degrees in $G$. We denote the common neighbourhood of nodes $i$ and $j$ by $S(i,j)=\Gamma_G(i)\cap\Gamma_G(j)$, so that $|S(i,j)|=(g^2)_{ij}$. For $i\in V$, let $G_i:=G[\Gamma_G(i)]$ denote the subgraph of $G$ induced by the neighbourhood of $i$, and let
$g^{(i)}$ denote the adjacency matrix of $G_i$. For $r\in V(G_i)=\Gamma_G(i)$, define $k_r^{(i)} := \sum_{s\in\Gamma_G(i)} g_{rs}$. Thus $k_r^{(i)}$ is the degree of node $r$ within the neighbourhood subgraph $G_i$, and need not equal its degree $k_r$ in $G$. For distinct nodes $i$ and $j$, let $G_{ij}:=G[S(i,j)]$ denote the subgraph induced by their common neighbourhood, and let $g^{(i,j)}$ denote its adjacency matrix. Equivalently, when
$j\in\Gamma_G(i)$, $G_{ij}$ is the subgraph of $G_i$ induced by the
neighbours of $j$ within $G_i$.\\
\\
\noindent
Let $G'=(V',E')\subseteq G$ denote a subgraph of $G$, with
$V'\subseteq V$ and $E'\subseteq E$. The subgraph $G'$ is induced if
$E'$ contains every edge $\{i,j\}\in E$ for which $i,j\in V'$.
Unless stated otherwise, our subgraph counts do not impose this
inducedness condition, so additional edges between the selected nodes
are permitted. For a graph pattern $H$ and a host graph $J$, let
$N(H,J)$ denote the number of general copies of $H$ in $J$. Special graphs used below include the complete graph on $n$ nodes, $K_n$,
which contains all possible edges, and the Erd{\H{o}}s--R{\'e}nyi
random graph $G(n,p)$ on node set $\{1,\ldots,n\}$, with edges that
arise independently with constant probability $p$.\\
\\
\noindent
We use the notation $M_a^{(b)}$ of \cite{agasse-duval_lawford25,lawford_mehmeti20} for a specific $b$-node subgraph pattern. To define $a$, label the nodes $1,\ldots,b$ and read the upper triangle of the adjacency matrix row by row as a $\binom{b}{2}$-bit binary number, with the first entry as the most significant bit. For a relabelling $\pi\in S_b$, let $h^\pi$ denote the resulting adjacency matrix. The canonical adjacency code is
\begin{equation*}
    a = \min_{\pi\in S_b}
    \sum_{i=1}^{b-1}\sum_{j=i+1}^{b}
    2^{\binom{b-i}{2}+(b-j)}(h^\pi)_{ij}.
\end{equation*}
Thus, for fixed $b$, $a$ uniquely identifies the isomorphism class of
the subgraph and is independent of its node labelling. Let $\widetilde{M}_a^{(b)}$ denote the corresponding induced pattern. For
the ambient graph $G$, we use the shorthand $|M_a^{(b)}| := N(M_a^{(b)},G)$ for the number of general copies of $M_a^{(b)}$, and write $|\widetilde{M}_a^{(b)}|$ for the number of induced copies of
$\widetilde{M}_a^{(b)}$ in $G$.

\begin{example}
Consider the spinning top:
\begin{center}
\includegraphics[scale=0.24]{M_119_5}
\end{center}

\noindent
Its $5!$ node relabellings produce the following 60 distinct decimal
adjacency codes:
\begin{equation*}
\begin{split}
\mathcal{A} ={}&
\{119,125,126,175,187,190,231,249,287,315,
317,343,378,399,444,461,\\
&462,467,470,473,483,485,490,500,543,559,567,605,622,
667,694,711,717,\\
&723,730,732,745,746,753,756,811,821,839,846,857,858,
867,873,876,882,\\
&884,903,918,921,924,933,938,940,945,946\}.
\end{split}
\end{equation*}
Hence $\min\mathcal{A}=119$, with
$119_{10}=0001110111_2$, and the spinning top is denoted
$M_{119}^{(5)}$. For later use, label its nodes
$i_1,\ldots,i_5$ so that $i_1$ has degree 1, $i_2$ has degree 2,
$i_5$ is the degree-3 node adjacent to $i_1$, and $i_3$ and $i_4$
are the remaining degree-3 nodes, in either order.
\end{example}

\section{The main result}\label{sec:main_result}
Our main result presents analytic count formulae for all five-node
subgraphs, in terms of functions of node degrees, powers of the
ambient adjacency matrix, adjacency matrices of neighbourhood and
common-neighbourhood subgraphs, and counts of smaller subgraphs. A
full combinatorial proof follows the statement of
Theorem~\ref{thm:analytic_count_5node}.

\begin{theorem}[Count formulae for subgraphs on five nodes]\label{thm:analytic_count_5node}

\begin{equation}\label{eq:m_75_5}
|M_{75}^{(5)}| = \sum_{i:k_i>3}\binom{k_i}{4} = \frac{1}{24}\sum_{i:k_i>3}k_i(k_i-1)(k_i-2)(k_i-3).
\end{equation}

\begin{equation}\label{eq:m_77_5}
\begin{split}
|M_{77}^{(5)}|
={}&\sum_{\{i,j\}\in E}\left[\binom{k_i-1}{2}(k_j-1)+
\binom{k_j-1}{2}(k_i-1)\right]-2|M_{15}^{(4)}| \\
={}&\frac{1}{2}\sum_{\{i,j\}\in E}\left[
(k_i-1)(k_i-2)(k_j-1)+(k_j-1)(k_j-2)(k_i-1)\right]
-2|M_{15}^{(4)}|.
\end{split}
\end{equation}

\begin{equation}\label{eq:m_79_5}
|M_{79}^{(5)}| = \frac{1}{2}\sum_{i:k_{i} > 3}(g^{3})_{ii} \dbinom{k_{i}-2}{2} = \frac{1}{4} \, \sum_{i:k_{i}>3}(g^{3})_{ii}(k_{i}-2)(k_{i}-3).
\end{equation}

\begin{equation}\label{eq:m_86_5}
|M_{86}^{(5)}| = \frac{1}{2}\sum_{i\neq j}(g^4)_{ij}-2|M_{3}^{(3)}|-9|M_{7}^{(3)}|-3|M_{11}^{(4)}|-2|M_{13}^{(4)}|-2|M_{15}^{(4)}|.
\end{equation}

\begin{equation}\label{eq:m_87_5}
|M_{87}^{(5)}|=\sum_{\substack{\{i,j\}\in E\\k_i>2,\,k_j>2}}
(g^2)_{ij}(k_i-2)(k_j-2)-2|M_{31}^{(4)}|.
\end{equation}

\begin{equation}\label{eq:m_94_5}
|M_{94}^{(5)}| = \sum_{\substack{i\neq j\\k_i>2}}\binom{(g^2)_{ij}}{2}(k_i-2)-2|M_{31}^{(4)}|=\frac{1}{2}\sum_{\substack{i\neq j\\k_i>2}}
(g^2)_{ij}\bigl((g^2)_{ij}-1\bigr)(k_i-2)-2|M_{31}^{(4)}|.
\end{equation}

\begin{equation}\label{eq:m_95_5}
|M_{95}^{(5)}| = \sum_{\substack{(i,j):\,\{i,j\}\in E\\k_i>3}}
\binom{(g^2)_{ij}}{2}(k_i-3) = \frac{1}{2} \sum_{\substack{(i,j):\,\{i,j\}\in E\\k_i>3}} (g^2)_{ij}\bigl((g^2)_{ij}-1\bigr)(k_i-3).
\end{equation}

\begin{equation}\label{eq:m_117_5}
|M_{117}^{(5)}| = \frac{1}{2}\sum_{(i,j):\,\{i,j\}\in E}(g^{3})_{ii}(k_j-1)-6|M_{7}^{(3)}|-2|M_{15}^{(4)}|-4|M_{31}^{(4)}|.
\end{equation}

\begin{equation}\label{eq:m_119_5}
|M_{119}^{(5)}| = \sum_{\substack{\{i,j\}\in E\\k_i>2,\,k_j>2}}
\bigl((g^2)_{ij}-1\bigr)\sum_{r\in S(i,j)}(k_r-2)-12|M_{63}^{(4)}|.
\end{equation}

\begin{equation}\label{eq:m_127_5}
|M_{127}^{(5)}| = \sum_{i:k_i>3}N(M_{7}^{(3)},G_i)\,(k_i-3)
= \frac{1}{6}\sum_{i:k_i>3}\tr\!\left((g^{(i)})^{3}\right)(k_i-3).
\end{equation}

\begin{equation}\label{eq:m_222_5}
|M_{222}^{(5)}|=\frac{1}{2}\sum_{\substack{i\neq j\\k_i>2,\,k_j>2}}
\binom{(g^2)_{ij}}{3}=\frac{1}{12}\sum_{\substack{i\neq j\\k_i>2,\,k_j>2}}(g^2)_{ij}\bigl((g^2)_{ij}-1\bigr)
\bigl((g^2)_{ij}-2\bigr).
\end{equation}

\begin{equation}\label{eq:m_223_5}
\left|M_{223}^{(5)}\right| = \sum_{\{i,j\}\in E}
\binom{(g^2)_{ij}}{3} = \frac{1}{6} \sum_{\{i,j\}\in E}
(g^2)_{ij} \bigl((g^2)_{ij}-1\bigr)\bigl((g^2)_{ij}-2\bigr).
\end{equation}

\begin{equation}\label{eq:m_235_5}
|M_{235}^{(5)}| = \sum_{i:k_{i}>2}\binom{\frac{1}{2}(g^{3})_{ii}}{2} - 2 \, |M_{31}^{(4)}| = \frac{1}{8} \, \sum_{i:k_{i}>2}(g^{3})_{ii}((g^{3})_{ii}-2) - 2 \, |M_{31}^{(4)}|.
\end{equation}

\begin{equation}\label{eq:m_236_5}
|M_{236}^{(5)}| = \frac{1}{10}\left(\tr(g^{5}) - 30 \, |M_{7}^{(3)}| - 10 \, |M_{15}^{(4)}|\right).
\end{equation}

\begin{equation}\label{eq:m_237_5}
|M_{237}^{(5)}| = \sum_{\{i,j\}\in E}(g^{3})_{ij}(g^{2})_{ij}-9|M_{7}^{(3)}|-2|M_{15}^{(4)}|-4|M_{31}^{(4)}|.
\end{equation}

\begin{equation}\label{eq:m_239_5}
|M_{239}^{(5)}| = \sum_{i:k_i>3}N(M_{13}^{(4)},G_i)
= \sum_{i:k_i>3}\left(\sum_{\{j,r\}\in E(G_i)}
(k_j^{(i)}-1)(k_r^{(i)}-1)-3N(M_{7}^{(3)},G_i)
\right).
\end{equation}

\begin{equation}\label{eq:m_254_5}
|M_{254}^{(5)}| = \sum_{\{i,j\}\in E} \;
\sum_{\{r,q\}\subset S(i,j)}\bigl((g^2)_{rq}-2\bigr).
\end{equation}

\begin{equation}\label{eq:m_255_5}
|M_{255}^{(5)}| = \frac{1}{2}\sum_{i:k_i>3}N(M_{15}^{(4)},G_i)
= \frac{1}{4}\sum_{i:k_i>3}\sum_{\substack{r\in V(G_i)\\k_r^{(i)}>2}}
\bigl((g^{(i)})^{3}\bigr)_{rr}\bigl(k_r^{(i)}-2\bigr).
\end{equation}

\begin{equation}\label{eq:m_507_5}
|M_{507}^{(5)}| = \sum_{i:k_i>3}N(M_{30}^{(4)},G_i)
=\frac{1}{8}\sum_{i:k_i>3}\left(\tr\!\left((g^{(i)})^{4}\right)
-2m(G_i)-4N(M_{3}^{(3)},G_i)\right).
\end{equation}

\begin{equation}\label{eq:m_511_5}
|M_{511}^{(5)}| = \frac{1}{3}\sum_{i:k_i>3}N(M_{31}^{(4)},G_i)
=\frac{1}{3}\sum_{i:k_i>3}\sum_{\{r,s\}\in E(G_i)}
\binom{((g^{(i)})^{2})_{rs}}{2}.
\end{equation}

\begin{equation}\label{eq:m_1023_5}
|M_{1023}^{(5)}| = \frac{1}{5}\sum_{i:k_i>3}N(M_{63}^{(4)},G_i)
=\frac{1}{120}\sum_{i:k_i>3}\sum_{\substack{j\in V(G_i)\\k_j^{(i)}>2}}
\tr\!\left((g^{(i,j)})^3\right).
\end{equation}

\end{theorem}

\noindent
Before proving the theorem, we record a general identity underlying the neighbourhood-subgraph construction used in several of the formulae. Recall that $N(H,J)$ denotes the number of general copies of a graph pattern $H$ in a host graph $J$. The graph $H$ need not be connected. Thus, for the subgraphs considered in
Theorem~\ref{thm:analytic_count_5node}, $N(M_a^{(b)},G)=|M_a^{(b)}|$.

\begin{lemma}[Neighbourhood lifting]\label{lem:neighbourhood_lifting}
Let $F$ be a finite simple graph and let $r\in V(F)$ be adjacent to
every other vertex of $F$. Set $H=F-r$, and let
\[
    \rho(F,r) := |\operatorname{Orb}_F(r)|,
\]
where $\operatorname{Orb}_F(r)$ denotes the orbit of $r$ under the
automorphism group of $F$. Then, for every simple graph $G$,
\[
    N(F,G) = \frac{1}{\rho(F,r)}\sum_{i\in V(G)} N(H,G_i).
\]
\end{lemma}

\begin{proof}
Let $\mathcal{C}_F(G)$ denote the set of copies of $F$ in $G$. For a
copy $C\in\mathcal{C}_F(G)$, let
\[
    R(C) = \bigl\{v\in V(C):
        \text{there exists an isomorphism }
        \varphi:F\to C
        \text{ with }\varphi(r)=v
    \bigr\}.
\]
For any fixed isomorphism $\varphi:F\to C$, the set $R(C)$ is the image
under $\varphi$ of $\operatorname{Orb}_F(r)$. Hence
\[
    |R(C)|=\rho(F,r).
\]
Now double-count the set
\[
    \mathcal{Q}
    =
    \{(C,i):C\in\mathcal{C}_F(G),\ i\in R(C)\}.
\]
Counting first by copies of $F$ gives
\[
    |\mathcal{Q}|=\rho(F,r)N(F,G).
\]
Alternatively, fix $i\in V(G)$. Since $r$ is adjacent to every other vertex of $F$, deleting $i$ from any pair $(C,i)\in\mathcal{Q}$ leaves a copy of $H=F-r$ whose vertices all lie in $\Gamma_G(i)$, and hence a copy of $H$ in $G_i$. Conversely, if $D$ is a copy of $H$ in $G_i$, then adjoining $i$ and the edges from $i$ to every vertex of $D$ produces a copy of $F$ in $G$ in which $i$ plays the role of $r$. These two operations are inverse to one another. Therefore the number of pairs $(C,i)\in\mathcal{Q}$ with second component $i$ is $N(H,G_i)$. Summing over $i\in V(G)$ gives
\[
    |\mathcal{Q}| = \sum_{i\in V(G)}N(H,G_i).
\]
Equating the two expressions for $|\mathcal{Q}|$ proves the result.
\end{proof}

\noindent The factor $\rho(F,r)$ records the number of structurally equivalent vertices of a copy of $F$ that can play the role of the
distinguished vertex $r$. Lemma~\ref{lem:neighbourhood_lifting} gives
alternative derivations of the crown, lamp, arrowhead, cat's cradle,
and 5-complete identities below, for which the corresponding values of
$\rho(F,r)$ are $1$, $2$, $1$, $3$, and $5$, respectively. We retain
the direct combinatorial proofs because they make the individual
constructions and multiplicity corrections explicit. We now prove Theorem~\ref{thm:analytic_count_5node}, grouping the identities according to the four recurring combinatorial constructions
described above.

\begin{proof}[Proof of Theorem \ref{thm:analytic_count_5node}]
The following identities are obtained using four recurring  constructions, based on incident structures, walks, common neighbourhoods, and neighbourhood subgraphs. See Table~\ref{tab:subgraphs34} for the connected subgraphs on three and four nodes.

\begin{table}[h]
\begin{center}
\begin{tabular}{cccc}
\includegraphics[height=0.8in]{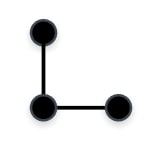} &
\includegraphics[height=0.8in]{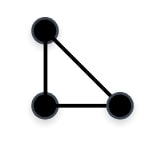} &
\includegraphics[height=0.8in]{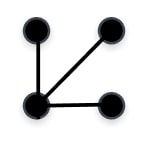} &
\includegraphics[height=0.8in]{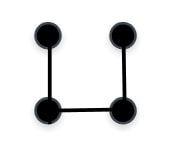} \\
3-star $M_{3}^{(3)}$ & triangle $M_{7}^{(3)}$ & 4-star $M_{11}^{(4)}$ & 4-path $M_{13}^{(4)}$\\
\\
\\
\includegraphics[height=0.8in]{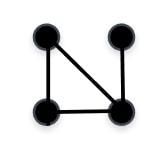} &
\includegraphics[height=0.8in]{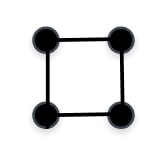} &
\includegraphics[height=0.8in]{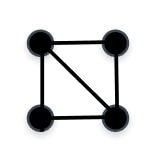} &
\includegraphics[height=0.8in]{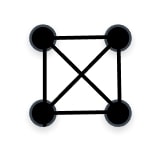} \\
tadpole $M_{15}^{(4)}$ & 4-circle $M_{30}^{(4)}$ & diamond $M_{31}^{(4)}$ & 4-complete $M_{63}^{(4)}$\\
\\
\end{tabular}
\caption{The 8 non-isomorphic connected subgraphs on three and four nodes.}
\label{tab:subgraphs34}
\end{center}
\end{table}

\medskip
\noindent\textbf{Incident structures.}
Begin with a distinguished node or edge and count smaller
structures incident to it, with corrections for multiple counting and
unwanted configurations.

\medskip
\noindent\emph{5-star $|M_{75}^{(5)}|$, equation~\eqref{eq:m_75_5}:\\\noindent}For each node $i$ with $k_i>3$, any four distinct neighbours of $i$ form a 5-star centred on $i$. Since the central degree-4 node is unique, each 5-star is counted once, giving \eqref{eq:m_75_5}. See also \cite[Proposition A.1, eqn.~18]{lawford_mehmeti20}.

\medskip
\noindent\emph{5-arrow $|M_{77}^{(5)}|$, equation~\eqref{eq:m_77_5}:\\\noindent}Consider an edge $\{i,j\}\in E$, and first let $i$ play the degree-3 role and $j$ the degree-2 role in a 5-arrow. Two neighbours of $i$ other than $j$ can be chosen in $\binom{k_i-1}{2}$ ways, and the remaining neighbour of $j$ can be
chosen in $k_j-1$ ways. Interchanging the roles of $i$ and $j$ gives
the second term in \eqref{eq:m_77_5}. The only unwanted case occurs
when the additional neighbour of $j$ coincides with one of the two
selected neighbours of $i$, producing a tadpole $M_{15}^{(4)}$.
Each tadpole is generated twice in this way, and hence
$2|M_{15}^{(4)}|$ is subtracted. This gives \eqref{eq:m_77_5}. See also \cite[Proposition A.1, eqn.~19]{lawford_mehmeti20}.

\medskip
\noindent\emph{cricket $|M_{79}^{(5)}|$,
equation~\eqref{eq:m_79_5}:\\\noindent}
The method of proof follows that used for the count of the tadpole
$|M_{15}^{(4)}|$ in \cite{agasse-duval_lawford25}. The cricket
subgraph can be thought of as a triangle on nodes $i$, $j$ and $x$,
with the addition of two extra edges $\{i,y\}$ and $\{i,z\}$, where
$k_i>3$. The quantity $(1/2)(g^3)_{ii}$ is the number of triangles
that contain node $i$, where the division by two corrects for the two
possible directions of travel around each triangle. Once a triangle
containing $i$ is fixed, two neighbours of $i$ have already been used,
and any two of its remaining $k_i-2$ neighbours can play the roles of
$y$ and $z$. Hence there are
\[
    \frac{1}{2}(g^3)_{ii}\binom{k_i-2}{2}
\]
crickets centred on node $i$. Summing over all eligible nodes $i$
gives \eqref{eq:m_79_5}. See also \cite[{$n_G(H_7)$}]{alon_etal97}.

\medskip
\noindent\emph{bull $|M_{87}^{(5)}|$, equation~\eqref{eq:m_87_5}:\\\noindent}The proof uses similar ideas to the count of the 4-path $|M_{13}^{(4)}|$ in \cite{agasse-duval_lawford25}. Consider an edge $\{i,j\}\in E$ as the central edge of a bull. Given this edge, $(g^2)_{ij}$ is the number of common neighbours of $i$ and $j$, and hence the number of triangles that contain $\{i,j\}$. Fix one such common neighbour $x$. Node $i$ then has $k_i-2$ possible remaining neighbours for node $y$, and node $j$ has $k_j-2$ possible remaining neighbours for node $z$, where $k_i>2$ and $k_j>2$. Hence $(g^2)_{ij}(k_i-2)(k_j-2)$ candidate bulls are associated with the central edge $\{i,j\}$. Summing over all eligible edges gives the first term in \eqref{eq:m_87_5}. This sum includes the unwanted case $y=z$, which forms a diamond with central edge $\{i,j\}$. Each diamond is counted twice, since either of its two common neighbours of $i$ and $j$ can play the role of $x$. We therefore subtract $2|M_{31}^{(4)}|$, giving \eqref{eq:m_87_5}. See also \cite[{$n_G(H_8)$}]{alon_etal97}.

\medskip
\noindent\emph{lollipop $|M_{117}^{(5)}|$, equation~\eqref{eq:m_117_5}:\\\noindent}The proof uses similar ideas to that of the count of the nested 5-arrow $|M_{77}^{(5)}|$ in \cite{lawford_mehmeti20}. Consider an edge $\{i,j\}\in E$, and let $i$ and $j$ play the degree-3 and degree-2 roles, respectively, in a lollipop. Let $i$ be directly connected to nodes $x$ and $z$, which form a triangle with $i$, and let $j$ also be directly connected to node $y$. The quantity $(1/2)(g^3)_{ii}$ is the number of triangles that contain node $i$, where the division by two corrects for the two possible directions of travel around each triangle. Node $j$ has $k_j-1$ possible remaining neighbours for node $y$. Hence there are $(1/2)(g^3)_{ii}(k_j-1)$ candidate lollipops associated with the orientation in which $i$ plays the degree-3 role and $j$ the degree-2 role. Summing over both orientations of every edge gives the first term in \eqref{eq:m_117_5}. This sum includes the unwanted cases $x=y$ and $z=y$, both of which form a diamond. Since four of the five edges of a diamond, with the appropriate orientation, can play the role of the central edge $\{i,j\}$ of a candidate lollipop, we subtract $4|M_{31}^{(4)}|$. Furthermore, each edge of a triangle
$M_{7}^{(3)}$, in both orientations, and two edges of a tadpole
$M_{15}^{(4)}$, with the appropriate orientation, can play the role
of the central edge of a candidate lollipop. We therefore also
subtract $6|M_{7}^{(3)}|$ and $2|M_{15}^{(4)}|$, giving
\eqref{eq:m_117_5}. See also \cite[{$n_G(H_{10})$}]{alon_etal97}.

\medskip
\noindent\emph{kite $|M_{127}^{(5)}|$, equation~\eqref{eq:m_127_5}:\\\noindent}Fix a node $i$ that plays the unique degree-4 role in the kite. Each triangle in the neighbourhood subgraph $G_i$ forms, together with $i$, a 4-complete subgraph containing $i$. There are $N(M_{7}^{(3)},G_i) = (1/6)\tr((g^{(i)})^{3})$ such triangles. For each one, any of the remaining $k_i-3$ neighbours of $i$ can play the degree-1 role in the kite. Hence there are $(1/6)\tr((g^{(i)})^{3})(k_i-3)$ kites whose unique degree-4 node is $i$. Summing over all eligible nodes $i$ gives \eqref{eq:m_127_5}.

\medskip
\noindent\emph{hourglass $|M_{235}^{(5)}|$,
equation~\eqref{eq:m_235_5}:\\\noindent}Let $i$ play the central degree-4 role in an hourglass. The quantity $(1/2)(g^3)_{ii}$ is the number of triangles that contain node $i$, where the division by two corrects for the two possible directions of travel around each triangle. Choosing any two such triangles gives $\binom{\frac{1}{2}(g^3)_{ii}}{2}$ candidate hourglasses associated with node $i$. Summing over all nodes with $k_i>2$ gives the total number of candidates. The unwanted case occurs when the two triangles share an edge, producing a diamond $M_{31}^{(4)}$. Each diamond is counted twice, once at each endpoint of its central edge. We therefore subtract $2|M_{31}^{(4)}|$, giving \eqref{eq:m_235_5}. See also
\cite[{$n_G(H_{11})$}]{alon_etal97}.

\medskip
\noindent\textbf{Walks.}
Count walks between one or more distinguished nodes,
sometimes together with additional incident structure, and correct
for unwanted cases.

\medskip
\noindent\emph{5-path $|M_{86}^{(5)}|$, equation~\eqref{eq:m_86_5}:\\\noindent}The quantity $(1/2)\sum_{i\neq j}(g^4)_{ij}$ counts walks of length four between distinct endpoints, with division by two correcting for the two possible orientations. Removing the walks that revisit nodes gives the correction terms in
\eqref{eq:m_86_5}. See \cite[Theorem 4.1]{movarraei_shikare14} and
\cite[Proposition A.1, eqn.~20]{lawford_mehmeti20}.

\medskip
\noindent\emph{banner $|M_{94}^{(5)}|$, equation~\eqref{eq:m_94_5}:\\\noindent}Let $i$ and $j$ play the degree-3 and degree-2 roles, respectively, in a banner. The quantity $(g^2)_{ij}$ gives the number of walks of length 2 between $i$ and $j$. Choosing any two distinct such walks gives $\binom{(g^2)_{ij}}{2}$ candidate 4-cycles through $i$ and $j$. A banner is then formed by adjoining to $i$ one of its remaining $k_i-2$ neighbours, say $x$. Hence
\[
    \binom{(g^2)_{ij}}{2}(k_i-2)
\]
candidate banners are associated with the ordered pair $(i,j)$. Summing over all ordered pairs $i\neq j$ gives the first term in
\eqref{eq:m_94_5}. The unwanted case $x=j$ produces a diamond. Each
diamond is generated twice, once for each endpoint of its central
edge playing the role of $i$, so we subtract $2|M_{31}^{(4)}|$.
This gives \eqref{eq:m_94_5}. See also \cite[{$n_G(H_9)$}]{alon_etal97}; the known typographical omission in that formula is discussed in
Example~\ref{example:further_comparison_with_published}.

\medskip
\noindent\emph{stingray $|M_{95}^{(5)}|$, equation~\eqref{eq:m_95_5}:\\\noindent}Let $i$ and $j$ play the degree-4 and degree-3 roles, respectively, in a stingray, with $\{i,j\}\in E$. Given this common edge, $(g^2)_{ij}$ is the number of common neighbours of $i$ and $j$, and hence the number of triangles containing $\{i,j\}$. Choosing any two of these common neighbours gives $\binom{(g^2)_{ij}}{2}$ pairs of triangles with common edge
$\{i,j\}$. A stingray is formed by adjoining to node $i$ one of its
remaining $k_i-3$ neighbours. Hence
\[
    \binom{(g^2)_{ij}}{2}(k_i-3)
\]
stingrays are associated with the orientation in which $i$ plays the
degree-4 role and $j$ the degree-3 role. Summing over both orientations of every eligible edge gives \eqref{eq:m_95_5}.

\medskip
\noindent\emph{ufo $|M_{222}^{(5)}|$,
equation~\eqref{eq:m_222_5}:\\\noindent}The method of proof is similar to that used for the count of the diamond $|M_{31}^{(4)}|$ in \cite{agasse-duval_lawford25}. We can think of a ufo on nodes $i$ and $j$ (both degree-3), and $x$, $y$ and $z$ (all degree-2), as three distinct paths of length two between $i$ and $j$. The quantity $(g^2)_{ij}$ is the number of common neighbours of $i$ and $j$, equivalently the number of walks of length 2 between them. A ufo is formed by choosing any three distinct common neighbours, so $\binom{(g^2)_{ij}}{3}$ gives the number of candidate ufos associated with the ordered pair $(i,j)$. Summing over all ordered pairs $i\neq j$ counts each ufo twice, once as $(i,j)$ and once as $(j,i)$. Dividing by two gives \eqref{eq:m_222_5}.

\medskip
\noindent\emph{chevron $|M_{223}^{(5)}|$,
equation~\eqref{eq:m_223_5}:\\\noindent}The method of proof is similar to that used for the count of the diamond $|M_{31}^{(4)}|$ in \cite{agasse-duval_lawford25}, and follows immediately from the proof of the ufo count $|M_{222}^{(5)}|$. We can think of a chevron on nodes $i$ and $j$ (both degree-4), and $x$, $y$ and $z$ (all degree-2), as three distinct triangles with a common edge $\{i,j\}$. Given this common edge, $(g^2)_{ij}$ is the number of walks of length 2 between $i$ and $j$, that is, the number of distinct triangles in $G$ that contain $\{i,j\}$. A chevron is formed by any three of these triangles, and so $\binom{(g^2)_{ij}}{3}$ gives the number of distinct chevrons that can be built from the common edge $\{i,j\}$. Summing over all edges $\{i,j\}\in E$ gives the total number of chevrons in $G$, which gives \eqref{eq:m_223_5}. See also
\cite[{$n_G(H_{13})$}]{alon_etal97}.

\medskip
\noindent\emph{5-circle $|M_{236}^{(5)}|$,
equation~\eqref{eq:m_236_5}:\\\noindent}The proof proceeds as for the 4-circle count $|M_{30}^{(4)}|$ in \cite{agasse-duval_lawford25}. The element $(g^5)_{ii}$ is the number of walks of length 5 from node $i$ back to itself, so $\tr(g^5)$ is the total number of closed walks of length 5 in $G$. A 5-circle contributes 10 such walks, corresponding to the five possible starting nodes and the two possible directions of travel. There are two further types of closed walks of length 5. First, each tadpole $M_{15}^{(4)}$ contributes 10 such walks: two from its degree-1 node, four from its degree-3 node, and two from each of its two degree-2 nodes. Second, each triangle $M_{7}^{(3)}$ contributes 30 such walks, five from each node in each of the two directions. Hence $\tr(g^5) = 10|M_{236}^{(5)}| + 10|M_{15}^{(4)}| + 30|M_{7}^{(3)}|$, and \eqref{eq:m_236_5} follows. See also
\cite[Theorem~2]{harary_manvel71} and \cite[{$n_G(C_5)$}]{alon_etal97}.

\medskip
\noindent\emph{house $|M_{237}^{(5)}|$,
equation~\eqref{eq:m_237_5}:\\\noindent}Let $\{i,j\}\in E$, and let $i$ and $j$ play the two degree-3 roles in a house. A house is formed by one walk of each of lengths 1, 2 and 3 between $i$ and $j$. For the fixed edge $\{i,j\}$, there are $(g^2)_{ij}$ walks of length 2 and $(g^3)_{ij}$ walks of length 3 between its endpoints. Hence $(g^3)_{ij}(g^2)_{ij}$ candidate houses are associated with the edge $\{i,j\}$. This count includes three unwanted cases. First, if $\{i,j\}$ is an edge of a triangle, there is one walk of length 2 and three walks of length 3 between $i$ and $j$, giving a total contribution of $9|M_{7}^{(3)}|$. Second, in a tadpole $M_{15}^{(4)}$, there are two edges joining a degree-2 role to the degree-3 role that generate an unwanted candidate, giving $2|M_{15}^{(4)}|$. Third, in a diamond $M_{31}^{(4)}$, there are four edges joining a degree-2 role to a degree-3 role that generate an unwanted candidate, giving $4|M_{31}^{(4)}|$. Subtracting these cases gives \eqref{eq:m_237_5}. See also \cite[{$n_G(H_{12})$}]{alon_etal97}.

\medskip
\noindent\textbf{Common neighbourhoods.}
Exploit the common neighbourhood
$S(i,j)=\Gamma_G(i)\cap\Gamma_G(j)$ of two distinguished nodes.

\medskip
\noindent\emph{spinning top $|M_{119}^{(5)}|$,
equation~\eqref{eq:m_119_5}:\\\noindent}Given an edge $\{i,j\}\in E$, let $S(i,j)$ be the common neighbourhood of nodes $i$ and $j$, with cardinality $|S(i,j)|=(g^2)_{ij}$. Choose a node $r\in S(i,j)$ and a second node $r'\in S(i,j)\setminus\{r\}$. The nodes $i,j,r,r'$ form a
diamond with central edge $\{i,j\}$. A candidate spinning top is
obtained by adjoining a further neighbour $x$ of $r$, distinct from
$i$ and $j$. There are $k_r-2$ choices for $x$, and $(g^2)_{ij}-1$ choices for $r'$. Hence the number of candidate spinning tops associated with the edge $\{i,j\}$ is $((g^2)_{ij}-1)\sum_{r\in S(i,j)}(k_r-2)$. The unwanted case occurs when $x=r'$, in which case the four nodes $i,j,r,r'$ form a 4-complete subgraph $M_{63}^{(4)}$. Each 4-complete subgraph generates 12 such unwanted candidates: there are six choices for its central edge $\{i,j\}$ and, for each one, two
choices for which of the remaining nodes plays the role of $r$. We therefore subtract $12|M_{63}^{(4)}|$, giving \eqref{eq:m_119_5}.

\medskip
\noindent\emph{envelope $|M_{254}^{(5)}|$,
equation~\eqref{eq:m_254_5}:\\\noindent}Let $\{i,j\}\in E$ and let $S(i,j)$ denote the common neighbourhood of nodes $i$ and $j$. Choose an unordered pair $\{r,q\}\subset S(i,j)$. The nodes $i,j,r,q$ form a diamond, with $\{i,j\}$ as its central edge. An envelope is obtained by adjoining a two-step path from $r$ to $q$ through another node $x$.
There are $(g^2)_{rq}$ common neighbours of $r$ and $q$. Two of these
are $i$ and $j$, so there are $(g^2)_{rq}-2$ possible choices for $x$. Summing over all unordered pairs $\{r,q\}\subset S(i,j)$ and
then over all edges $\{i,j\}\in E$ gives \eqref{eq:m_254_5}.

\medskip
\noindent\textbf{Neighbourhood subgraphs.}
Count smaller subgraph patterns inside the neighbourhood
subgraph $G_i$ and correct for multiple counting.

\medskip
\noindent\emph{crown $|M_{239}^{(5)}|$,
equation~\eqref{eq:m_239_5}:\\\noindent}
Fix a node $i$ that plays the unique degree-4 role in the crown. Any
4-path $M_{13}^{(4)}$ contained in the neighbourhood subgraph $G_i$
forms, together with $i$, a crown $M_{239}^{(5)}$. Hence
$N(M_{13}^{(4)},G_i)$ is the number of crowns whose unique degree-4
node is $i$, and summing over all eligible nodes $i$ gives the first
equality in \eqref{eq:m_239_5}. To obtain the second equality, apply
the 4-path formula
\[
|M_{13}^{(4)}| = \sum_{\{j,r\}\in E}(k_j-1)(k_r-1) - 3|M_{7}^{(3)}|
\]
from \cite{agasse-duval_lawford25} to the neighbourhood subgraph
$G_i$. The relevant degrees are therefore the local degrees
$k_j^{(i)}$ and $k_r^{(i)}$, rather than the corresponding degrees
in $G$. This gives
\[
N(M_{13}^{(4)},G_i) = \sum_{\{j,r\}\in E(G_i)}
(k_j^{(i)}-1)(k_r^{(i)}-1) - 3N(M_{7}^{(3)},G_i),
\]
and substitution gives the second equality in
\eqref{eq:m_239_5}. The first equality also follows directly from
Lemma~\ref{lem:neighbourhood_lifting}. Taking $F=M_{239}^{(5)}$ and
letting $r$ be its unique degree-4 node gives
$F-r\cong M_{13}^{(4)}$ and $\rho(F,r)=1$. Hence
\[
|M_{239}^{(5)}| = \sum_{i\in V(G)}N(M_{13}^{(4)},G_i).
\]

\medskip
\noindent\emph{lamp $|M_{255}^{(5)}|$,
equation~\eqref{eq:m_255_5}:\\\noindent}
Fix a node $i$ that plays one of the two structurally equivalent degree-4 roles in the lamp. Any tadpole $M_{15}^{(4)}$ contained in the neighbourhood subgraph $G_i$ forms, together with $i$, a lamp $M_{255}^{(5)}$. Hence $N(M_{15}^{(4)},G_i)$ counts the lamps in which $i$ plays one of the two degree-4 roles. Summing over all eligible nodes $i$ therefore counts each lamp twice, and dividing by two gives the first equality in \eqref{eq:m_255_5}. To obtain the second equality, apply the tadpole formula
\[
|M_{15}^{(4)}| = \frac{1}{2}\sum_{r:k_r>2}(g^3)_{rr}(k_r-2)
\]
from \cite{agasse-duval_lawford25} to the neighbourhood subgraph $G_i$. The relevant degree is then the local degree $k_r^{(i)}$, and the adjacency matrix is $g^{(i)}$. Thus
\[
N(M_{15}^{(4)},G_i) = \frac{1}{2} \sum_{\substack{r\in V(G_i)\\k_r^{(i)}>2}} \bigl((g^{(i)})^3\bigr)_{rr}
\bigl(k_r^{(i)}-2\bigr),
\]
which gives the second equality in \eqref{eq:m_255_5}. Equivalently, Lemma~\ref{lem:neighbourhood_lifting} gives the first equality directly. For $F=M_{255}^{(5)}$, either of its two structurally equivalent degree-4 nodes can serve as the distinguished vertex, and deleting either one leaves a tadpole $M_{15}^{(4)}$. Hence $\rho(F,r)=2$ and
\[
|M_{255}^{(5)}| = \frac{1}{2}\sum_{i\in V(G)}N(M_{15}^{(4)},G_i).
\]

\medskip
\noindent\emph{arrowhead $|M_{507}^{(5)}|$,
equation~\eqref{eq:m_507_5}:\\\noindent}
Fix a node $i$ that plays the unique degree-4 role in the arrowhead. Any 4-circle $M_{30}^{(4)}$ contained in the neighbourhood subgraph $G_i$ forms, together with $i$, an arrowhead $M_{507}^{(5)}$. Hence $N(M_{30}^{(4)},G_i)$ is the number of arrowheads whose unique degree-4 node is $i$, and summing over all eligible nodes $i$ gives the first equality in \eqref{eq:m_507_5}. To obtain the second equality, apply the 4-circle formula
\[
|M_{30}^{(4)}| = \frac{1}{8}\left(\tr(g^4)-2m-4|M_{3}^{(3)}|\right)
\]
from \cite{agasse-duval_lawford25} to the neighbourhood subgraph
$G_i$. Its adjacency matrix is $g^{(i)}$ and it has $m(G_i)$ edges,
which gives the second equality in \eqref{eq:m_507_5}. The first equality also follows directly from Lemma~\ref{lem:neighbourhood_lifting}. Taking $F=M_{507}^{(5)}$ and letting $r$ be its unique degree-4 node gives $F-r\cong M_{30}^{(4)}$ and $\rho(F,r)=1$. Therefore
\[
|M_{507}^{(5)}| = \sum_{i\in V(G)}N(M_{30}^{(4)},G_i).
\]

\medskip
\noindent\emph{cat's cradle $|M_{511}^{(5)}|$,
equation~\eqref{eq:m_511_5}:\\\noindent}Fix a node $i$ that plays one of the three structurally equivalent degree-4 roles in the cat's cradle. Deleting $i$ from the cat's cradle leaves a diamond $M_{31}^{(4)}$. Hence every diamond contained in the neighbourhood subgraph $G_i$ forms, together with $i$, a cat's cradle. Since each cat's cradle has three structurally equivalent degree-4 nodes, summing over $i$ counts each copy three times. Therefore
\[
|M_{511}^{(5)}| = \frac{1}{3}\sum_{i:k_i>3}N(M_{31}^{(4)},G_i).
\]
For a fixed $i$, a diamond in $G_i$ is determined by its unique central edge $\{r,s\}\in E(G_i)$ and two distinct common neighbours
of $r$ and $s$ in $G_i$. Their number is therefore
\[
N(M_{31}^{(4)},G_i) = \sum_{\{r,s\}\in E(G_i)}
\binom{((g^{(i)})^{2})_{rs}}{2}.
\]
Substitution gives the second equality in \eqref{eq:m_511_5}. Equivalently, the first equality follows directly from
Lemma~\ref{lem:neighbourhood_lifting}, with $F=M_{511}^{(5)}$, $F-r\cong M_{31}^{(4)}$, and $\rho(F,r)=3$.

\medskip
\noindent\emph{5-complete $|M_{1023}^{(5)}|$,
equation~\eqref{eq:m_1023_5}:\\\noindent}
Fix a node $i$ that plays one of the five structurally equivalent degree-4 roles in the 5-complete graph. Any 4-complete subgraph $M_{63}^{(4)}$ contained in the neighbourhood subgraph $G_i$ forms, together with $i$, a 5-complete graph. Since each 5-complete graph has five structurally equivalent vertices,
\[
|M_{1023}^{(5)}| = \frac{1}{5}\sum_{i:k_i>3}N(M_{63}^{(4)},G_i).
\]
To obtain the second equality, apply the 4-complete formula
\[
|M_{63}^{(4)}| = \frac{1}{24}\sum_{j:k_j>2} \tr\!\left((g^{(j)})^3\right)
\]
from \cite{agasse-duval_lawford25} to the neighbourhood subgraph $G_i$. The relevant degree of $j$ is therefore its local degree
$k_j^{(i)}$, and the relevant neighbourhood adjacency matrix is
$g^{(i,j)}$. Thus
\[
N(M_{63}^{(4)},G_i) = \frac{1}{24}
\sum_{\substack{j\in V(G_i)\\k_j^{(i)}>2}}
\tr\!\left((g^{(i,j)})^3\right),
\]
which gives the second equality in \eqref{eq:m_1023_5}. See also
\cite[Proposition A.1, eqn. 21]{lawford_mehmeti20}. Equivalently, the
first equality follows directly from
Lemma~\ref{lem:neighbourhood_lifting}, with $F=M_{1023}^{(5)}$,
$F-r\cong M_{63}^{(4)}$, and $\rho(F,r)=5$. 

\end{proof}

\begin{remark}[Results, proof strategies and extensions]\label{rem:results_strategies_extensions}
We make some observations on the techniques used to derive
Theorem~\ref{thm:analytic_count_5node}, on alternative approaches to
counting particular subgraphs, and on extensions of these methods to
larger subgraphs.
\begin{enumerate}[label=(\roman*)]
    \item\label{4} The proofs of equations~\eqref{eq:m_75_5}--\eqref{eq:m_1023_5} use four recurring strategies. These strategies are not mutually exclusive, and several identities admit more than one derivation:
    \begin{enumerate}
        \item {[Incident structures]} Start from a node $i$ or an edge $\{i,j\}\in E$ and choose one or more structures incident to it, such as nodes, triangles, or a 4-complete subgraph. This strategy is used for the 5-star $M_{75}^{(5)}$, 5-arrow $M_{77}^{(5)}$, cricket $M_{79}^{(5)}$, bull $M_{87}^{(5)}$, lollipop $M_{117}^{(5)}$, kite $M_{127}^{(5)}$, and hourglass $M_{235}^{(5)}$.
        
        \item {[Walks]} Start from nodes $i$ and $j$, not necessarily
        directly connected, and count walks between them, possibly together with additional structure incident to one or both nodes. This strategy is used for the 5-path $M_{86}^{(5)}$, banner $M_{94}^{(5)}$, stingray $M_{95}^{(5)}$, ufo $M_{222}^{(5)}$, chevron $M_{223}^{(5)}$, house $M_{237}^{(5)}$, and, when $i=j$, the 5-circle $M_{236}^{(5)}$.
        
        \item {[Common neighbourhoods]} Start from nodes $i$ and $j$ and their common neighbourhood $S(i,j)$, and then count structures involving members of $S(i,j)$. This strategy is used for the spinning top $M_{119}^{(5)}$ and envelope $M_{254}^{(5)}$.
        
        \item\label{5} {[Neighbourhood subgraphs]}
        Apply Lemma~\ref{lem:neighbourhood_lifting} by choosing a node $r$ of the target graph $F$ that is adjacent to every other node, deleting $r$, and counting the resulting graph $H=F-r$ inside $G_i$. The factor $\rho(F,r)^{-1}$ corrects for structurally equivalent choices of the deleted node. This strategy is used explicitly for the crown $M_{239}^{(5)}$ (4-path), lamp $M_{255}^{(5)}$ (tadpole), arrowhead $M_{507}^{(5)}$ (4-circle), cat's cradle $M_{511}^{(5)}$ (diamond), and 5-complete $M_{1023}^{(5)}$ (4-complete subgraph).
    \end{enumerate}
    
    In each case, we sum over the relevant nodes or edges, correcting for
multiple counting and removing ``unwanted cases'' that satisfy the same
local restrictions as the subgraph of interest. For example, a
candidate 5-arrow $M_{77}^{(5)}$ consists of an edge $\{i,j\}\in E$
together with two neighbours of $i$ and one neighbour of $j$; a
tadpole $M_{15}^{(4)}$ satisfies the same restrictions and must
therefore be removed twice. Table~\ref{tab:unwanted_cases} records the
smaller subgraphs and unwanted cases that appear explicitly in
equations~\eqref{eq:m_75_5}--\eqref{eq:m_1023_5}.\\
\\
\noindent The proof strategies employed by \cite{alon_etal97,estrada11,estrada_knight15} are closely related. These works cover the same eight five-node subgraphs: the cricket, bull, banner, lollipop, chevron, hourglass, house, and 5-circle, while the 5-star count is either given explicitly or follows immediately from the 4-star count. See also Example~\ref{example:further_comparison_with_published} for discussion of some of these formulae. Their count formulae can be derived using the incident-structure strategy for the 5-star, cricket, bull, lollipop, and hourglass, and the walk strategy for the banner, chevron, house, and 5-circle. The common-neighbourhood and neighbourhood-subgraph strategies are not used there in the systematic form considered here. Our use of the latter was motivated by the observation in \cite[p.~222]{alon_etal97} that 4-complete subgraphs can be counted by looking for triangles in the subgraph induced by the neighbourhood of a node. Here we extend this idea by counting smaller subgraph patterns inside neighbourhood subgraphs $G_i$.\\
\\
The formulae for the cricket, bull, banner, lollipop, chevron,
hourglass, house, and 5-circle appear in earlier work, primarily
\cite{alon_etal97}; the 5-star count is either given explicitly or
follows immediately from the corresponding 4-star count. We proved the
5-arrow, 5-path, and 5-complete formulae in \cite{lawford_mehmeti20}. For the remaining nine subgraphs (the stingray, spinning top, kite, ufo, crown, envelope, lamp, arrowhead, and cat's cradle) the analytical formulae given here appear to provide substantially new formulations. Under the grouping above, these nine results use all four strategies: incident structures for the kite; walks for the stingray and ufo; common neighbourhoods for the spinning top and envelope; and neighbourhood subgraphs for the crown, lamp, arrowhead, and cat's cradle. The 5-complete formula is a direct generalization of the observation in \cite[p.~222]{alon_etal97} that 4-complete subgraphs can be counted through triangles in neighbourhood
subgraphs. Thus each of the four strategies yields at least one new or
substantially reformulated five-node identity.
    
\begin{table}[p]
\centering
\scriptsize

\renewcommand{\arraystretch}{1.65}

\(
\begin{blockarray}{rcccccccc}
& |M_{3}^{(3)}| & |M_{7}^{(3)}| & |M_{11}^{(4)}| & |M_{13}^{(4)}|
& |M_{15}^{(4)}| & |M_{30}^{(4)}| & |M_{31}^{(4)}| & |M_{63}^{(4)}|\\
\begin{block}{r(cccccccc)}
  \textrm{5-star} \; |M_{75}^{(5)}| & \cdot & \cdot & \cdot & \cdot & \cdot & \cdot & \cdot & \cdot \\
  \textrm{5-arrow} \; |M_{77}^{(5)}| & \cdot & \cdot & \cdot & \cdot & \hlight{X} & \cdot & \cdot & \cdot \\
  \textrm{cricket} \; |M_{79}^{(5)}| & \cdot & \cdot & \cdot & \cdot & \cdot & \cdot & \cdot & \cdot \\
  \textrm{5-path} \; |M_{86}^{(5)}| & \hlight{X} & \hlight{X} & \hlight{X} & \hlight{X} & \hlight{X} & \cdot & \cdot & \cdot \\
  \textrm{bull} \; |M_{87}^{(5)}| & \cdot & \cdot & \cdot & \cdot & \cdot & \cdot & \hlight{X} & \cdot \\
  \textrm{banner} \; |M_{94}^{(5)}| & \cdot & \cdot & \cdot & \cdot & \cdot & \cdot & \hlight{X} & \cdot \\
  \textrm{stingray} \; |M_{95}^{(5)}| & \cdot & \cdot & \cdot & \cdot & \cdot & \cdot & \cdot & \cdot \\
  \textrm{lollipop} \; |M_{117}^{(5)}| & \cdot & \hlight{X} & \cdot & \cdot & \hlight{X} & \cdot & \hlight{X} & \cdot \\
  \textrm{spinning top} \; |M_{119}^{(5)}| & \cdot & \cdot & \cdot & \cdot & \cdot & \cdot & \cdot & \hlight{X} \\
  \textrm{kite} \; |M_{127}^{(5)}| & \cdot & \hlight{X} & \cdot & \cdot & \cdot & \cdot & \cdot & \cdot \\
  \textrm{ufo} \; |M_{222}^{(5)}| & \cdot & \cdot & \cdot & \cdot & \cdot & \cdot & \cdot & \cdot \\
  \textrm{chevron} \; |M_{223}^{(5)}| & \cdot & \cdot & \cdot & \cdot & \cdot & \cdot & \cdot & \cdot \\
  \textrm{hourglass} \; |M_{235}^{(5)}| & \cdot & \cdot & \cdot & \cdot & \cdot & \cdot & \hlight{X} & \cdot \\
  \textrm{5-circle} \; |M_{236}^{(5)}| & \cdot & \hlight{X} & \cdot & \cdot & \hlight{X} & \cdot & \cdot & \cdot \\
  \textrm{house} \; |M_{237}^{(5)}| & \cdot & \hlight{X} & \cdot & \cdot & \hlight{X} & \cdot & \hlight{X} & \cdot \\
  \textrm{crown} \; |M_{239}^{(5)}| & \cdot & \cdot & \cdot & \hlight{X} & \cdot & \cdot & \cdot & \cdot \\
  \textrm{envelope} \; |M_{254}^{(5)}| & \cdot & \cdot & \cdot & \cdot & \cdot & \cdot & \cdot & \cdot \\
  \textrm{lamp} \; |M_{255}^{(5)}| & \cdot & \cdot & \cdot & \cdot & \hlight{X} & \cdot & \cdot & \cdot \\
  \textrm{arrowhead} \; |M_{507}^{(5)}| & \cdot & \cdot & \cdot & \cdot & \cdot & \hlight{X} & \cdot & \cdot \\
  \textrm{cat's cradle} \; |M_{511}^{(5)}| & \cdot & \cdot & \cdot & \cdot & \cdot & \cdot & \hlight{X} & \cdot \\
  \textrm{5-complete} \; |M_{1023}^{(5)}| & \cdot & \cdot & \cdot & \cdot & \cdot & \cdot & \cdot & \hlight{X} \\
\end{block}
\end{blockarray}
\)

\caption{Smaller subgraphs appearing as direct components of the count
formulae (\ref{eq:m_75_5})--(\ref{eq:m_1023_5}). Subgraphs that arise
only through expansion of these component counts are not shown. The
five-node subgraphs and the three- and four-node subgraphs are shown in Tables \ref{tab:subgraphs5} and \ref{tab:subgraphs34}, respectively.}
\label{tab:unwanted_cases}

\end{table}

    \item Observe that Lemma~\ref{lem:neighbourhood_lifting} also gives an alternative formulation for the chevron $M_{223}^{(5)}$. Either of its two degree-4 nodes can serve as the distinguished vertex $r$, and deleting either one leaves a 4-star $M_{11}^{(4)}$. Hence $\rho(M_{223}^{(5)},r)=2$ and
    \[
        |M_{223}^{(5)}| = \frac{1}{2}\sum_{i\in V(G)} N(M_{11}^{(4)},G_i).
    \]
    Since the number of 4-stars in $G_i$ is
    \[
        N(M_{11}^{(4)},G_i) = \sum_{s\in V(G_i)}\binom{k_s^{(i)}}{3},
    \]
    we obtain
    \begin{equation}\label{eq:m_223_5_neighbourhood}
    |M_{223}^{(5)}| = \frac{1}{2}\sum_{i:k_i>3}\sum_{s\in V(G_i)}
    \binom{k_s^{(i)}}{3} = \frac{1}{12}\sum_{i:k_i>3}\sum_{s\in V(G_i)}k_s^{(i)}\bigl(k_s^{(i)}-1\bigr)\bigl(k_s^{(i)}-2\bigr).
    \end{equation}

    We have found some preliminary evidence in simulations that, for
    moderately sized sparse Erd{\H{o}}s--R{\'e}nyi graphs, the chevron can be counted faster using equation~\eqref{eq:m_223_5_neighbourhood}
    than using equation~\eqref{eq:m_223_5}.

    \item Lemma~\ref{lem:neighbourhood_lifting} also applies when the graph remaining after deletion of $r$ is disconnected. For example, let
    $F=M_{127}^{(5)}$ be the kite and let $r$ be its unique degree-4 node. Then
    \[
        F-r\cong K_3\sqcup K_1,
    \]
    where $\sqcup$ denotes disjoint union and $K_1$ denotes the graph
    consisting of a single node. Since $\rho(F,r)=1$,
    Lemma~\ref{lem:neighbourhood_lifting} gives
    \[
        |M_{127}^{(5)}| = \sum_{i\in V(G)}N(K_3\sqcup K_1,G_i).
    \]
    For each triangle in $G_i$, any one of the remaining $k_i-3$ nodes of $G_i$ can play the isolated role in $K_3\sqcup K_1$, because no inducedness condition is imposed. Therefore
    \[
        |M_{127}^{(5)}| = \sum_{i:k_i>3}N(M_7^{(3)},G_i)(k_i-3)
        = \frac{1}{6}\sum_{i:k_i>3}
        \operatorname{tr}\!\left((g^{(i)})^3\right)(k_i-3),
    \]
    which recovers \eqref{eq:m_127_5}.

    \item It is important to apply the strategies in item~\ref{4} carefully to avoid complications arising from unwanted cases. Consider the spinning top $M_{119}^{(5)}$ and let $\{i,j\}\in E$ be the edge between the two degree-3 nodes that are not directly connected to the degree-1 node. It might seem reasonable to exploit walks of lengths 1, 2 and 4 to count candidate spinning tops using
    \begin{equation*}
        \frac{1}{2}\sum_{i\neq j}
        \binom{(g^2)_{ij}g_{ij}}{2}(g^4)_{ij} =
        \sum_{\{i,j\}\in E}\binom{(g^2)_{ij}}{2}(g^4)_{ij}.
    \end{equation*}
    However, the unwanted cases now include five-node subgraphs
    (stingray, crown, chevron, and envelope), as well as larger subgraphs on six or seven nodes, so the argument no longer closes in terms of counts of simpler patterns.
    
    \item\label{rem:six_node_examples} The strategies in item~\ref{4} can be used to derive exact count formulae for larger subgraphs, and we give two examples.
    
    \begin{figure}\centering
        \begin{subfigure}{0.35\textwidth}
            \centering
            \includegraphics[scale = 0.3]{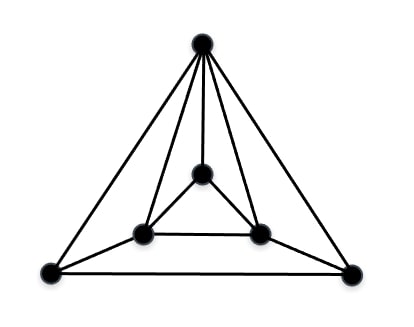}
            \caption{Subgraph $M_{7919}^{(6)}$.}
            \label{fig:M_7919_6}
        \end{subfigure}
        \begin{subfigure}{0.35\textwidth}
            \centering
            \includegraphics[scale = 0.4]{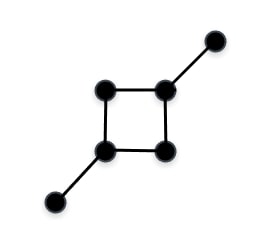}
            \caption{Subgraph $M_{1182}^{(6)}$.}
            \label{fig:M_1182_6}
        \end{subfigure}
        \caption{Two six-node subgraphs used to illustrate extensions of
        the proof strategies for Theorem~\ref{thm:analytic_count_5node}.}
        \label{fig:six_node_graphlets}
    \end{figure}
    
    \begin{enumerate}
    
    \item The neighbourhood-lifting construction is not restricted to
    five-node subgraphs. For example, consider the six-node subgraph
    $M_{7919}^{(6)}$ in Figure~\ref{fig:M_7919_6}. Its unique degree-5
    node may be taken as $r$, and deleting this node leaves a house
    $M_{237}^{(5)}$. Since $\rho(M_{7919}^{(6)},r)=1$,
    Lemma~\ref{lem:neighbourhood_lifting} gives
    \[
        |M_{7919}^{(6)}| = \sum_{i:k_i>4}N(M_{237}^{(5)},G_i).
    \]
    Thus the same identity applies directly to subgraphs with more
    than five nodes whenever the required structural condition is
    satisfied.
    
    \item Extending ideas from the banner $M_{94}^{(5)}$, we can count
    the six-node subgraph in Figure~\ref{fig:M_1182_6} using
    \begin{equation*}
        |M_{1182}^{(6)}| = \frac{1}{2}
        \sum_{\substack{i\neq j\\k_i>2,\,k_j>2}}
        \binom{(g^2)_{ij}}{2}(k_i-2)(k_j-2)
        -|M_{31}^{(4)}| -|M_{95}^{(5)}| -3|M_{222}^{(5)}|.
    \end{equation*}
    \end{enumerate}

    These examples illustrate the scope, but also the limitation, of the methods used here. Individual constructions often extend directly to larger subgraphs, and the neighbourhood-lifting identity applies whenever the target pattern contains an appropriate universal node. In general, however, the correction terms generated by the other constructions need not involve only smaller patterns, so the methods do not provide a closed recursive scheme for enumerating all connected subgraphs of a given larger order.
    
    \item Several authors present exact enumeration formulae for subgraphs on five nodes, or strategies for deriving combinatorial results. Partial results on five-node subgraphs, with proof, include \cite{alon_etal97} (8 subgraphs), \cite{harary_manvel71} (1 subgraph), and \cite{movarraei_shikare14} (1 subgraph). The eight formulae for five-node subgraphs listed, with proof strategies, in \cite[\S 4.4]{estrada11} and \cite[Chapter 13]{estrada_knight15}, and three of the four five-node subgraphs (excluding the 5-path) given without proof in \cite{allen-perkins_etal17}, all appear in \cite{alon_etal97}, essentially in the same form (although \cite{allen-perkins_etal17, estrada11} use matrix notation). To our knowledge, the only other complete set of results for five-node subgraphs is given by \cite{pinar_etal17}, although they use a very different method of proof and, in some cases, obtain quite different formulations. We discuss some related issues in Examples~\ref{example:alternative_forms} and
    \ref{example:further_comparison_with_published}, including the 5-path count given in \cite{allen-perkins_etal17}.
    \end{enumerate}
\end{remark}

\paragraph{Computational validation.}
As an independent computational check, we compared the 21 general
count formulae \eqref{eq:m_75_5}--\eqref{eq:m_1023_5} and the 21
induced count formulae in Appendix~A with direct brute-force
enumeration. The enumeration oracle implements the definitions of
general and induced subgraph copies directly and does not use any of
the analytical count identities. The comparison was exhaustive over
all $2^{15}=32{,}768$ labelled simple graphs on six nodes. We also
tested $60{,}000$ randomly generated Erd\H{o}s--R\'enyi graphs
$G(n,p)$, with $n=7,\ldots,12$, $p\in\{0.1,0.3,0.5,0.7,0.9\}$, and $2{,}000$ realizations for each $(n,p)$. All $92{,}768$ host graphs gave exact agreement between the analytical formulae and the independent enumeration.

\section{Examples}\label{sec:examples}
In this section, we present a series of illustrative theoretical and
empirical applications of the main result. First, we show how
Theorem~\ref{thm:analytic_count_5node} can be used to derive
non-trivial analytical results by specializing the count formulae for
the 5-path, bull, and spinning top to three regular graphs that have
been widely studied: the complete graph $K_n$, the complete
$N$-partite graph $K_{n_1,n_2,\ldots,n_N}$, and the regular ring
lattice (also known in the complex systems literature as the
small-world model with no edge rewiring \cite{watts_strogatz98}).
Second, we consider three different formulations of the 5-path count
formula from the literature, together with three results from
\cite{estrada11, estrada_knight15}, and show how easily
misunderstandings can arise. Third, we apply all 21 count formulae to
a real-world network, obtain the corresponding induced counts, and
compare the resulting induced five-node subgraph profile with
Erd\H{o}s--R\'enyi and degree-preserving null ensembles. We start by
proving a straightforward but useful lemma that gives the number of
walks of length $k$ between any pair of nodes in a complete graph
$K_n$, as a function of $k$ and $n$.

\begin{lemma}[Number of walks of length $k$ in a complete graph]
\label{thm:paths_k_complete}
Let $G=K_n$ be a complete graph with adjacency matrix $g$ and $n\geq 2$. For $k\geq 0$, write
\[
    a_k=(g^k)_{ii}, \quad b_k=(g^k)_{ij}, \quad i\neq j.
\]
Then
\begin{equation*}
    b_k=\frac{(n-1)^k+(-1)^{k+1}}{n}, \quad a_k=b_k+(-1)^k.
\end{equation*}
\end{lemma}

\begin{proof}
For $k\geq 3$, the off-diagonal elements of $g^k$ satisfy the
recursion
\[
    b_k=(n-2)b_{k-1}+(n-1)b_{k-2}, \quad b_1=1, \quad b_2=n-2.
\]
Indeed, matrix multiplication gives $b_k=a_{k-1}+(n-2)b_{k-1}$, while
$a_{k-1}=(n-1)b_{k-2}$. The characteristic polynomial of the
recursion is
\[
    \chi(x) = x^2-(n-2)x-(n-1) = (x+1)(x-n+1),
\]
with distinct roots $-1$ and $n-1$. Hence the general solution is
\[
    b_k = C(-1)^k+D(n-1)^k,
\]
where $C$ and $D$ are constants. Applying the initial conditions gives
\[
    (n-1)D-C=1, \quad (n-1)^2D + C =n-2,
\]
whose unique solution is
\[
    C=-1/n, \quad D=1/n,
\]
and the claim on $b_{k}$ follows immediately. This formula also holds for $k=0$, since the off-diagonal elements of $g^0=I$ are zero. For $k\geq 1$, we have $a_k=(n-1)b_{k-1}$, and we can easily verify that $a_{k} - b_{k} = (-1)^k$, and the claim on $a_k$ follows. The same identity holds for $k=0$, since $a_0=1$ and $b_0=0$.
\end{proof}

\begin{figure}\centering
    \begin{subfigure}[t]{0.3\textwidth}\vspace{0pt}
        \centering
        \includegraphics[width=1.0\linewidth]{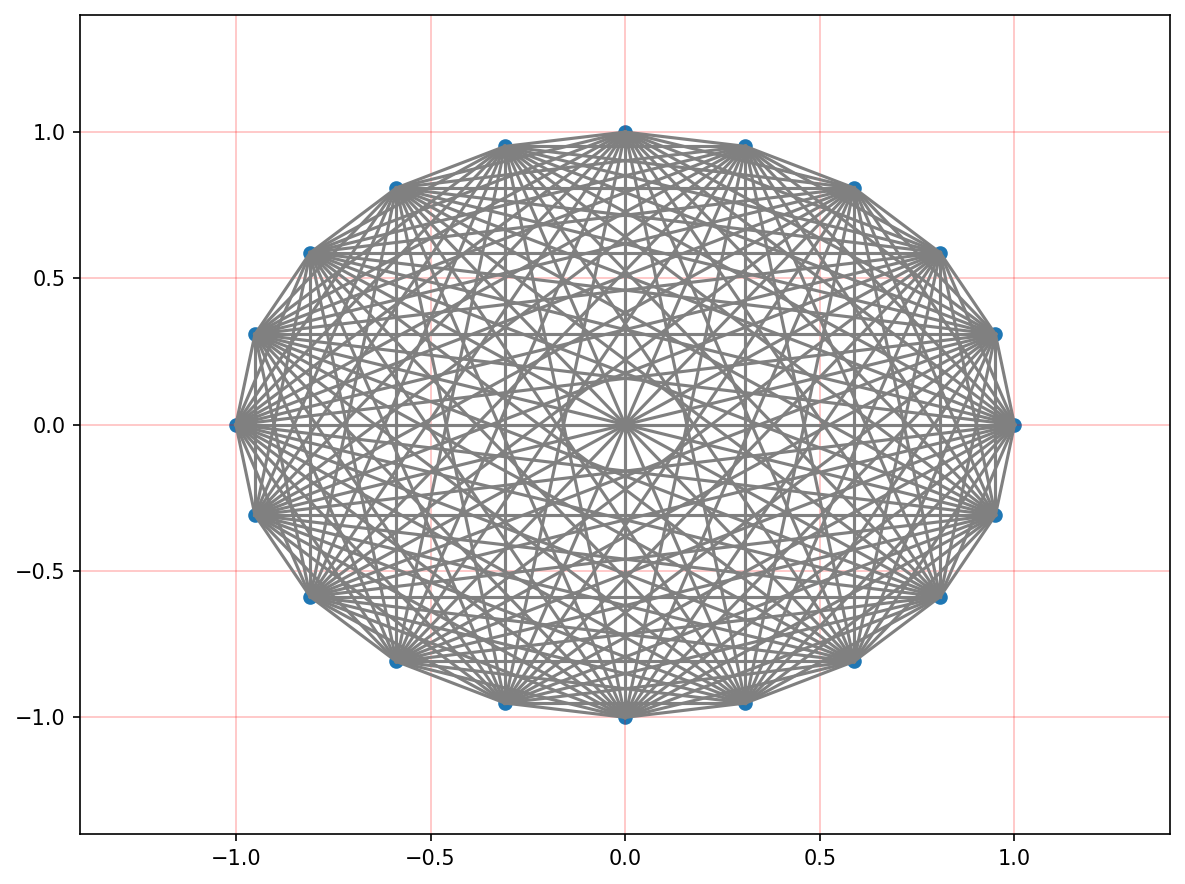}
        \caption{\scriptsize{Complete graph $K_{20}$.}}
        \label{fig:complete_K_20}
    \end{subfigure}
    \begin{subfigure}[t]{0.3\textwidth}\vspace{0pt}
        \centering
        \includegraphics[width=1.0\linewidth]{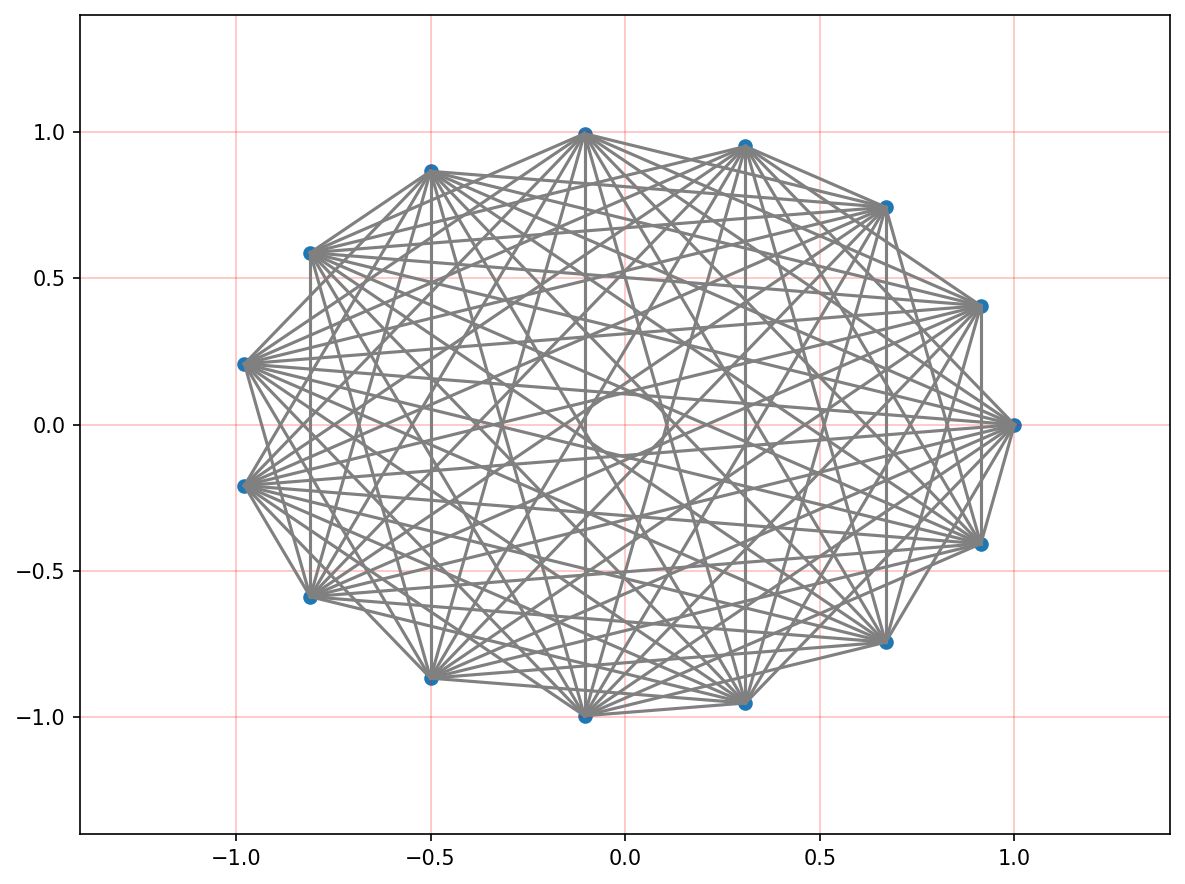}
        \caption{\scriptsize{Complete 5-partite graph.}}
        \label{fig:complete_5_partite}
    \end{subfigure}
    \begin{subfigure}[t]{0.3\textwidth}\vspace{0pt}
        \centering
        \includegraphics[width=1.0\linewidth]{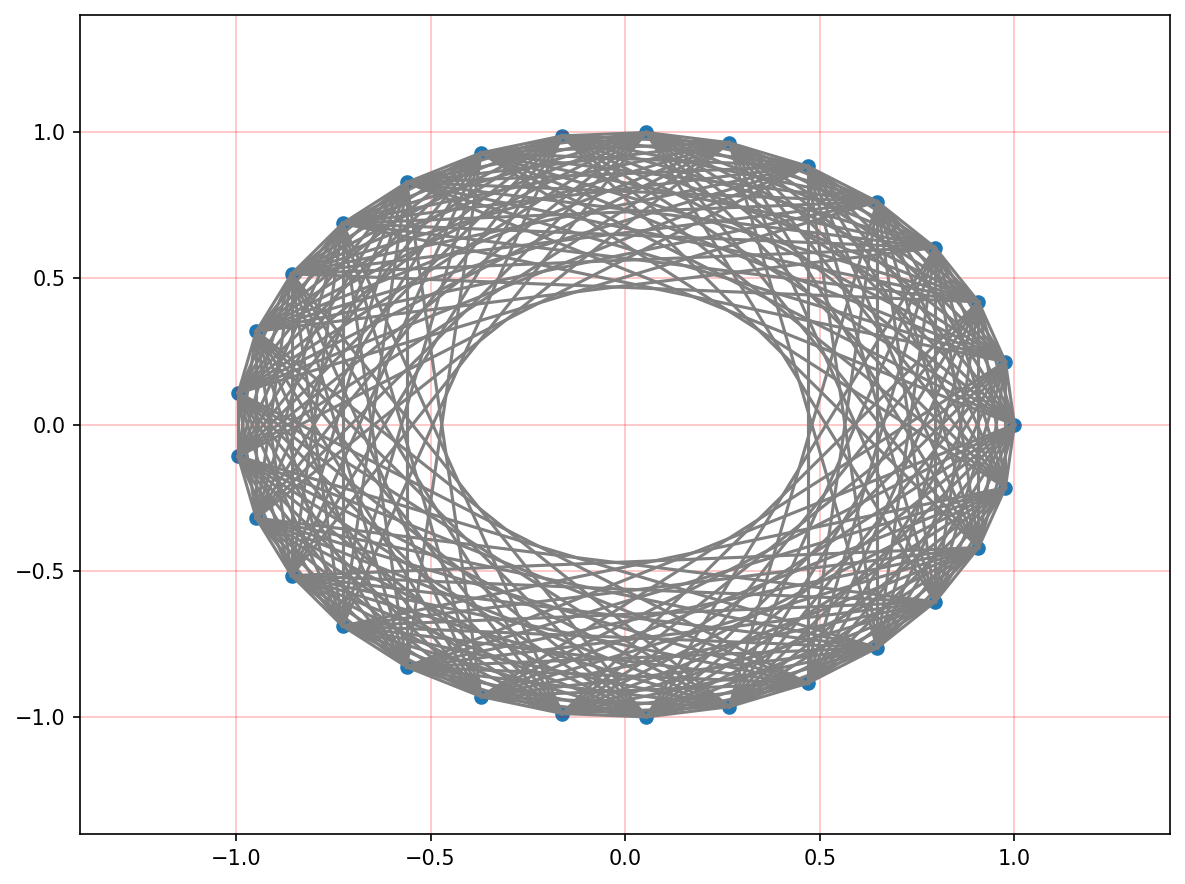}
        \caption{\scriptsize{Regular ring lattice.}}
        \label{fig:regular_ring_lattice}
    \end{subfigure}
    \caption{Three regular graphs used to illustrate analytic subgraph
    count formulae in Examples~\ref{example:5-path}--\ref{example:spinning_tops}.}
    \label{fig:illustrative_graphs}
    \end{figure}

\noindent As a sample application of Lemma~\ref{thm:paths_k_complete},
we specialize it to walks of length 4 in $K_n$. Directly, we have
$(g^4)_{ii}=(n-1)(n^2-3n+3)$ and $(g^4)_{ij}=(n-2)(n^2-2n+2)$. For example, in $K_{20}$, illustrated in Figure~\ref{fig:complete_K_20}, there are 6,517 walks of length 4 from any node to itself, and 6,516 walks of length 4 from any node to any other node.

\begin{example}[Counting 5-paths in the complete graph $K_n$]\label{example:5-path}
Assume $n\geq 5$. We compute each term of the 5-path count formula
\begin{equation*}
    |M_{86}^{(5)}| = \frac{1}{2}\sum_{i\neq j}(g^4)_{ij} - 2|M_{3}^{(3)}| - 9|M_{7}^{(3)}| - 3|M_{11}^{(4)}| - 2|M_{13}^{(4)}| - 2|M_{15}^{(4)}|
\end{equation*}
given in \eqref{eq:m_86_5}, noting that each node $i$ has degree
$k_i=n-1$. From Lemma~\ref{thm:paths_k_complete}, we can write the
first term as
\[
    \frac{1}{2}\sum_{i\neq j}(g^4)_{ij} = (n-2)(n^2-2n+2)\binom{n}{2}.
\]
The number of 3-stars follows from \cite[eqn.~1]{agasse-duval_lawford25} as
\[
    |M_{3}^{(3)}| = \sum_i\binom{k_i}{2} = n\binom{n-1}{2}.
\]
Applying Lemma~\ref{thm:paths_k_complete} and
\cite[eqn.~2]{agasse-duval_lawford25}, the number of triangles is
\[
    |M_{7}^{(3)}| = \frac{1}{6}\sum_i(g^3)_{ii} = \binom{n}{3}.
\]
From \cite[eqn.~3]{agasse-duval_lawford25}, the number of 4-stars is
\[
    |M_{11}^{(4)}| = \sum_i\binom{k_i}{3} = n\binom{n-1}{3}.
\]
Combining \cite[eqn.~4]{agasse-duval_lawford25} with the triangle
count above, the number of 4-paths is
\[
    |M_{13}^{(4)}| = \sum_{\{i,j\}\in E}(k_i-1)(k_j-1) - 3|M_{7}^{(3)}| = (n-2)^2\binom{n}{2} - 3\binom{n}{3}.
\]
Applying the recurrence
\[
    \binom{n}{k} = \frac{n-k+1}{k}\binom{n}{k-1}
\]
gives
\[
    |M_{13}^{(4)}| = (n-2)(n-3)\binom{n}{2}.
\]
Applying Lemma~\ref{thm:paths_k_complete} and
\cite[eqn.~5]{agasse-duval_lawford25} gives
\[
    |M_{15}^{(4)}| = \frac{1}{2}\sum_i(g^3)_{ii}(k_i-2) = 12\binom{n}{4}
\]
tadpoles. Using the binomial recurrence above twice gives
\[
    |M_{15}^{(4)}| = (n-2)(n-3)\binom{n}{2},
\]
and hence $|M_{13}^{(4)}|=|M_{15}^{(4)}|$ in $K_n$: there is the same
number of 4-paths and tadpoles in the complete graph. Substituting
these expressions into the 5-path formula and repeatedly applying the
binomial recurrences above gives
\[
    |M_{86}^{(5)}|=60\binom{n}{5}.
\]
The combinatorial intuition is immediate: there are $\binom{n}{5}$
ways to choose five nodes, and each set of five nodes supports
$5!/2=60$ distinct undirected 5-paths in $K_n$.
\end{example}

\begin{example}[Counting bulls in the complete $N$-partite graph]
\label{example:bull}
Let $G$ be a complete $N$-partite graph $K_{n_1,n_2,\ldots,n_N}$, with nodes partitioned into $N\geq 2$ disjoint sets (or \emph{groups}) such that no pair of nodes within the same group is adjacent, while all pairs of nodes in different groups are adjacent. Let
\[
    n:=\sum_{c=1}^{N}n_c,
\]
where $n_c$ is the number of nodes in group $c$. The graph $G$ has an
$n\times n$ block-symmetric adjacency matrix $g$, with block $(a,b)$
equal to $0_{n_a\times n_a}$ if $a=b$, and equal to
\[
    J_{a,b}:=\iota_{n_a}\iota_{n_b}'
\]
if $a\neq b$, for $a,b=1,2,\ldots,N$, where $\iota_q$ denotes the
$q$-vector of ones. Since
\[
    J_{a,b}J_{b,c}=n_bJ_{a,c},
\]
every element of block $(a,b)$ of $g^k$ has the same value. It is
therefore convenient to work with blocks rather than individual
nodes. Define
\[
    A^{(1)}:=\iota_N\iota_N'-I_N, \quad D:=\diag(n_1,n_2,\ldots,n_N).
\]
It is straightforward to verify that each element of block $(a,b)$ of
$g^k$ equals $(A^{(k)})_{ab}$, where
\[
    A^{(k)} = A^{(1)}DA^{(k-1)} = (A^{(1)}D)^{k-1}A^{(1)}.
\]
We compute the bull subgraph count formula
\begin{equation*}
|M_{87}^{(5)}| = \sum_{\substack{\{i,j\}\in E\\k_i>2,\,k_j>2}}
(g^2)_{ij}(k_i-2)(k_j-2) - 2|M_{31}^{(4)}|
\end{equation*}
given in \eqref{eq:m_87_5}, and use the diamond count
\[
|M_{31}^{(4)}| = \frac{1}{2} \sum_{\{i,j\}\in E}
(g^2)_{ij}\bigl((g^2)_{ij}-1\bigr)
\]
from \cite[eqn.~7]{agasse-duval_lawford25}. To simplify the example,
we assume that every node in group $c$ has degree $n-n_c>2$, so that every edge in the graph is potentially the ``central'' edge of a bull subgraph. Observe that
\[
    (A^{(2)})_{ab}=n-n_a-n_b \quad (a\neq b).
\]
The desired count may then be rewritten in terms of blocks as
\begin{equation*}
    |M_{87}^{(5)}| = \sum_{1\leq a<b\leq N} n_an_b(n-n_a-n_b)
    \left[(n-n_a-2)(n-n_b-2) - (n-n_a-n_b-1)\right].
\end{equation*}
To go further, suppose that all groups have the same size,
$n_1=\cdots=n_N=:s$, so that $n=Ns$. Then
\begin{equation}\label{eq:m_87_5_complete_N_partite}
    |M_{87}^{(5)}| = 3\binom{N}{3}s^3
    \left\{\bigl((N-1)s-2\bigr)\bigl((N-1)s-3\bigr)+s-1\right\}.
\end{equation}
For fixed $s$, $|M_{87}^{(5)}|=O(N^5)$, while for fixed $N$,
$|M_{87}^{(5)}|=O(s^5)$. Thus the number of bull subgraphs grows on
the order of the fifth power of either the number of groups or the
common group size when the other is held fixed. Trivially,
$|M_{87}^{(5)}|=0$ when $N\leq 2$ or $n<5$. For $N\geq 3$, the
bracketed factor in \eqref{eq:m_87_5_complete_N_partite} has
positive-integer zeros only at $(N,s)=(3,1)$ and $(4,1)$, both of
which have $n<5$. Hence every equal-part-size complete $N$-partite
graph with $N\geq 3$ and $n\geq 5$ contains at least one bull
subgraph.\\
\\
\noindent
We now consider the combinatorial intuition behind
\eqref{eq:m_87_5_complete_N_partite}. Choose three groups from the
$N$ groups, each containing one node of a triangle. There are $s^3$
possible triangles for any such choice of groups. We then choose one
of the three pairs of nodes in the triangle to play the two degree-3
roles in the bull; denote these nodes by $i$ and $j$. Thus there are
\[
    3\binom{N}{3}s^3
\]
ways to choose the triangle and the two degree-3 roles. Let $A$ denote the group containing $j$. There are $(N-1)s-2$ choices for a node $k$ adjacent to $i$ to form the first horn, after excluding the other two nodes of the triangle. If $k\notin A$, then $k$ is also adjacent to $j$, so there are $(N-1)s-3$ choices for a distinct node $\ell$ adjacent to $j$ to form the second horn. If $k\in A$, then $k$ is not adjacent to $j$ and therefore does not reduce the number of available neighbours of $j$; in this case there is one additional choice for $\ell$. Since $s-1$ of the possible choices of $k$ lie in $A$, the total number of ways to choose the two horns, for a fixed triangle and pair of degree-3 roles, is
\[
    \bigl((N-1)s-2\bigr)\bigl((N-1)s-3\bigr)+(s-1),
\]
which gives \eqref{eq:m_87_5_complete_N_partite}. For example, in the
complete 5-partite graph with $s=3$ nodes in each group, illustrated
in Figure~\ref{fig:complete_5_partite}, there are 74,520 bull
subgraphs.
\end{example}

\begin{example}[Counting spinning tops in the regular ring lattice]
\label{example:spinning_tops}
Let $G$ be a regular ring lattice on $n$ nodes, where each node is
adjacent to its $k$ nearest neighbours in both the clockwise and
anti-clockwise directions, giving degree $2k$ at every node and a
total of $nk$ edges. We require $n>2k$. We compute the spinning-top
subgraph formula
\begin{equation*}
|M_{119}^{(5)}| = \sum_{\substack{\{i,j\}\in E\\k_i>2,\,k_j>2}}
\bigl((g^2)_{ij}-1\bigr) \sum_{r\in S(i,j)}(k_r-2) - 12|M_{63}^{(4)}|
\end{equation*}
given in \eqref{eq:m_119_5}, where $|M_{63}^{(4)}|$ is the
4-complete subgraph count and $S(i,j)$ is the common neighbourhood
of nodes $i$ and $j$. We identify two cases corresponding to a
discrete change in the behaviour of the spinning-top count at the
threshold $n=3k+1$.\\
\\
\noindent
Case A: $n\geq 3k+1$. We start by considering $n$ sufficiently large
relative to $k$ that it is enough to consider the structure among a
node's $k$ nearest neighbours in the clockwise direction. Each node
forms a 4-complete subgraph with any three of these neighbours, and
each such subgraph is counted once in this way. Hence
\[
    |M_{63}^{(4)}|=n\binom{k}{3}.
\]
Every node $r$ has global degree $k_r=2k$. Since $g$ and $g^2$ are
both symmetric circulant matrices, it is enough to consider, for
$a=1,\ldots,k$, an edge $\{i,j\}$ for which $j$ is the $a$-th
nearest neighbour of $i$ in the clockwise direction. In this case,
\[
    |S(i,j)|=(g^2)_{ij}=2k-a-1.
\]
As $i$ varies, these choices give each edge of $G$ exactly once.
Substitution into \eqref{eq:m_119_5} gives
\begin{equation*}
|M_{119}^{(5)}| = 2(k-1)n \sum_{a=1}^{k} [2k-(a+1)][2k-(a+2)] -
12n\binom{k}{3},
\end{equation*}
which, after simplification, yields
\begin{equation}\label{eq:spinning_top_count_case_A}
|M_{119}^{(5)}| = \frac{2}{3}nk(k-1)\left[7(k-1)(k-2)+3\right].
\end{equation}
Thus the expression is linear in $n$ and quartic in $k$, or
$O(nk^4)$. Since the degree of each node is $2k$, the dependence on
the degree is also quartic.\\
\\
\noindent
Case B: $2k+1\leq n\leq 3k$. As $n$ decreases from $3k$ to $2k+1$,
there are progressively more 4-complete subgraphs incident to a given
node $i$, more paths of length 2 between the endpoints of each edge,
and more members of each common neighbourhood $S(i,j)$. Define
\[
    \alpha_a := \max\{a-(n-2k-1),0\},
\]
and
\[
    B:=n\,\mathrm{Te}_{3k-n}, \quad
    \mathrm{Te}_{\ell}:=\binom{\ell+2}{3},
\]
where $\mathrm{Te}_{\ell}$ is the $\ell$-th tetrahedral number. One
can verify that
\begin{equation*}
|M_{119}^{(5)}| = 2(k-1)n \sum_{a=1}^{k}[2k-(a+2)+\alpha_a]
[2k-(a+1)+\alpha_a] - 12\left[n\binom{k}{3}+B\right].
\end{equation*}
Separating out the Case A count in
\eqref{eq:spinning_top_count_case_A} gives
\begin{equation*}
|M_{119}^{(5)}| = \frac{2}{3}nk(k-1)[7(k-1)(k-2)+3] +
2(k-1)n\sum_{a=1}^{k}\left[\alpha_a(4k-2a-3)+\alpha_a^2\right]
- 12B.
\end{equation*}
Since $\alpha_a=0$ for $a\leq n-2k-1$, this becomes
\begin{align*}
    |M_{119}^{(5)}|
    ={}&
    \frac{2}{3}nk(k-1)
    [7(k-1)(k-2)+3] \\
    &+
    2(k-1)n
    \sum_{a=n-2k}^{k}
    \Bigl[
    (a-(n-2k-1))(4k-2a-3) \\
    &\hspace{5.5cm}
    +(a-(n-2k-1))^2
    \Bigr] \\
    &-
    12n\binom{3k-n+2}{3} \\[0.5ex]
    ={}&
    \frac{2}{3}n
    \Bigl\{
    (3k-n+1)(3k-n+2) \\
    &\hspace{1.6cm}\times
    [9k^2-(2n+21)k+5n+3] \\
    &\hspace{1.6cm}
    +[7(k-1)(k-2)+3](k-1)k
    \Bigr\}.
\end{align*}

\noindent Let $f(n)$ denote the subgraph count $|M_{119}^{(5)}|$ for
$n\geq2k+1$, conditional on $k$. When $k=1$, there are no spinning
tops for any $n$, so consider $k\geq2$. At the boundary between Cases
A and B, we have
\[
    f(3k+1)-f(3k) = \frac{2}{3}k(k-1)(7k^2-39k+35).
\]
This quantity is negative for $k=2,3,4$, so $f(n)$ is not monotonically increasing in $n$ for these values of $k$. For $k\geq5$,
\[
\begin{split}
    f(2k+1)-f(3k)
    &=
    2k(k-1)(k^3+9k^2-25k+9) \\
    &=
    2k(k-1)\bigl[k(k^2+9k-25)+9\bigr]
    >0.
\end{split}
\]
Since $2k+1<3k$, this again rules out monotonic increase. Thus
$|M_{119}^{(5)}|$ is not monotonically increasing in $n$ for any
$k\geq2$. Moreover, \eqref{eq:spinning_top_count_case_A} is linear
in $n$ with positive slope for $k\geq2$, so the number of spinning
tops is minimized for some $n$ satisfying $2k+1\leq n\leq3k+1$. For example, when $k=10$,
\[
    |M_{119}^{(5)}| =
    \begin{cases}
        488724n-39026n^2+1092n^3-10n^4,
        & 21\leq n\leq30,\\
        30420n,
        & n\geq31,
    \end{cases}
\]
and the number of spinning tops is minimized at $n=29$, with a count
of 912,108. This graph is illustrated in
Figure~\ref{fig:regular_ring_lattice}.
\end{example}

\begin{example}[Alternative formulations of the 5-path count formula]
\label{example:alternative_forms}
We now examine three different formulations of the 5-path count formula that we found in the applied mathematics literature. For reference, our result \eqref{eq:m_86_5} is
\begin{equation*}
    |M_{86}^{(5)}| = \frac{1}{2}\sum_{i\neq j}(g^4)_{ij} - 2|M_{3}^{(3)}| - 9|M_{7}^{(3)}| - 3|M_{11}^{(4)}|
    - 2|M_{13}^{(4)}| - 2|M_{15}^{(4)}|.
\end{equation*}
In each case, we begin by quoting the authors' result in its original
notation, without defining the terms, before restating it using the
notation of our paper. This comparison illustrates the importance of
clear and consistent notation for subgraphs.\\
\noindent \textbf{Formulation 1} (correct). We start with
\cite[Theorem 4.1]{movarraei_shikare14}, who formulate their result in
terms of individual matrix elements. Writing the left-hand side of
their expression as $A$, their result is
\begin{equation*}
    \text{``$A = \sum_{i \neq j}[a_{ij}^{(4)} - 2 \, a_{ij}^{(2)} \, (d_{j} - a_{ij})] - \sum_{i=1}^{n}[(2 \, d_{i} - 1) \, a_{ii}^{(3)} + 6 \, \binom{d_{i}}{3}]$''},
\end{equation*}
where $A$ denotes twice the 5-path count. Restating the right-hand
side in the notation of our paper gives
\begin{align*}
    \frac{A}{2}
    &=
    \frac{1}{2}
    \left(
        \sum_{i\neq j}
        \left[(g^4)_{ij}
        -2(g^2)_{ij}(k_j-g_{ij})\right]
        -
        \sum_i
        \left[(2k_i-1)(g^3)_{ii}
        +6\binom{k_i}{3}\right]
    \right) \\
    &=
    \frac{1}{2}\sum_{i\neq j}(g^4)_{ij}
    -
    \sum_{i\neq j}(g^2)_{ij}(k_j-g_{ij})
    -
    \sum_i k_i(g^3)_{ii}
    +
    \frac{1}{2}\sum_i(g^3)_{ii}
    -
    3\sum_i\binom{k_i}{3}.
\end{align*}
We use the following results. From \cite[eqn.~3]{agasse-duval_lawford25},
\[
    -3\sum_i\binom{k_i}{3} = -3|M_{11}^{(4)}|.
\]
From \cite[eqn.~2]{agasse-duval_lawford25},
\[
    \frac{1}{2}\sum_i(g^3)_{ii} = \frac{1}{2}\tr(g^3) = 3|M_{7}^{(3)}|.
\]
Finally, from \cite[eqn.~5]{agasse-duval_lawford25},
\[
    -\sum_i k_i(g^3)_{ii} = -2|M_{15}^{(4)}| - 12|M_{7}^{(3)}|.
\]
It therefore remains to show that
\begin{equation*}
    \sum_{i\neq j}(g^2)_{ij}(k_j-g_{ij}) =
    2\left(|M_{3}^{(3)}|+|M_{13}^{(4)}|\right).
\end{equation*}
Observe that
\begin{align*}
    \sum_{i\neq j}(g^2)_{ij}(k_j-g_{ij})
    &=
    \sum_{i\neq j}(g^2)_{ij} + \sum_{i\neq j}(g^2)_{ij}(k_j-g_{ij}-1).
\end{align*}
Combining this decomposition with \cite[\S 2.1 and Theorem 3.7]{movarraei_shikare14} establishes the identity above, and hence the equivalence of \eqref{eq:m_86_5} and \cite[Theorem 4.1]{movarraei_shikare14}.\\
\\
\noindent \textbf{Formulation 2} (as stated, this is correct only in a
special case). We now examine \cite[Lemma 9, eqn.~4.12]{allen-perkins_etal17}, who present their result as a linear combination of small subgraph counts,
\begin{equation*}
    \text{``$|P_{4}| = \frac{1}{2} \, \vec{1}^{T}A^{4}\vec{1} - |P_{1}| - 4 \, |P_{2}| - 2 \, |P_{3}| - 9 \, |C_{3}| - 4 \, |C_{4}| - 6 \, |S_{1,3}| - 4 \, |S_{T1S}|$''}.
\end{equation*}
In the notation of our paper, this becomes
\begin{equation*}
    |P_{4}| = \frac{1}{2}\sum_{i,j}(g^4)_{ij} - m - 4|M_{3}^{(3)}|
    - 9|M_{7}^{(3)}| - 6|M_{11}^{(4)}| - 2|M_{13}^{(4)}|
    - 4|M_{15}^{(4)}| - 4|M_{30}^{(4)}|.
\end{equation*}

\noindent Using \cite[eqn.~6]{agasse-duval_lawford25} to separate the diagonal terms in the matrix sum, this reduces to
\begin{equation*}
    |P_{4}| = \frac{1}{2}\sum_{i\neq j}(g^4)_{ij} - 2|M_{3}^{(3)}|
    - 9|M_{7}^{(3)}| - 6|M_{11}^{(4)}| - 2|M_{13}^{(4)}|
    - 4|M_{15}^{(4)}|.
\end{equation*}
Comparison with \eqref{eq:m_86_5} therefore shows that
$|P_{4}|=|M_{86}^{(5)}|$ if and only if
\begin{equation*}
    6|M_{11}^{(4)}|+4|M_{15}^{(4)}| = 3|M_{11}^{(4)}|+2|M_{15}^{(4)}|.
\end{equation*}
Since subgraph counts are non-negative, this holds if and only if
\[
    |M_{11}^{(4)}|=|M_{15}^{(4)}|=0.
\]
Equivalently, $G$ has maximum degree at most two. Thus, if $G$ is
connected, it is a path or a cycle; more generally, it is a disjoint
union of paths, cycles, and isolated nodes. The comparison suggests
that the coefficients of $|S_{1,3}|$ and $|S_{T1S}|$ in the stated
expression for $|P_4|$ have inadvertently been doubled.\\
\\
\noindent \textbf{Formulation 3} (as stated, this is incorrect even for a path). We end with \cite[Theorem 9, $N_3$]{pinar_etal17}, who also give a result in terms of smaller subgraph counts,
\begin{equation*}
    \text{``$N_{3} = \sum_{i}\sum_{(i,j)\in E}(d(j)-1)-4\cdot C_{4}(G)-2\cdot TT(G)-3\cdot T(G)$''}.
\end{equation*}
The double summation means that each edge is considered in both
directions. Rewriting the result in the notation of our paper gives
\begin{equation*}
    N_{3} = \sum_{\{i,j\}\in E}(k_i+k_j-2) - 3|M_{7}^{(3)}|
    - 2|M_{15}^{(4)}| - 4|M_{30}^{(4)}|.
\end{equation*}

\noindent We show by counterexample that this is false. Let $G=P_n$ be the path on $n\geq5$ nodes. Immediately,
\[
    |M_{7}^{(3)}| = |M_{11}^{(4)}| = |M_{15}^{(4)}| = 0,
    \quad |M_{3}^{(3)}|=n-2, \quad |M_{13}^{(4)}|=n-3.
\]
Hence, from \eqref{eq:m_86_5},
\begin{equation*}
    |M_{86}^{(5)}| = \frac{1}{2}\sum_{i\neq j}(g^4)_{ij} -4n + 10.
\end{equation*}
The sum $\sum_{i\neq j}(g^4)_{ij}$ counts non-closed walks of length
4 with their starting and ending nodes distinguished. In $P_n$, such
a walk can end only two or four nodes away from its starting node.
For pairs of nodes at distance two, there are three walks when the
pair is adjacent to an endpoint of the path and four walks otherwise.
For pairs at distance four, there is one walk in each direction.
Therefore
\begin{align*}
    \sum_{i\neq j}(g^4)_{ij}
    &=
    2\left[3+3+4(n-4)\right]+2(n-4) \\
    &=
    10n-28.
\end{align*}
It follows that
\[
    |M_{86}^{(5)}| = \frac{1}{2}(10n-28)-4n+10 = n-4,
\]
as expected. For the expression of \cite[Theorem 9, $N_3$]{pinar_etal17}, we have $|M_{30}^{(4)}|=0$, and hence
\[
    N_3 = \sum_{\{i,j\}\in E}(k_i+k_j-2) = 2n-4.
\]
Thus \eqref{eq:m_86_5} and \cite[Theorem 9, $N_3$]{pinar_etal17} are not equivalent on $P_n$ for any $n\geq5$. For example, when $G=P_5$ is a 5-path, $|M_{86}^{(5)}|=1$, as required, whereas $N_3=6$.
\end{example}

\begin{example}[Further comparison with published results]
\label{example:further_comparison_with_published}
It is sometimes difficult to verify published exact subgraph count
formulae, or to distinguish typographical errors from more substantive
errors. Some excellent papers explicitly omit proofs of the count
formulae altogether. For example,
\cite[p.~267]{allen-perkins_etal17} state that ``The proofs of these
results are based on the strategy developed and explained in
\cite{estrada_knight15} [our reference] and are not given here as they
are lengthy and technical.'' Similarly, \cite[p.~1436]{pinar_etal17}
write that ``It is cumbersome and space-consuming to give proofs of all
of these, so we omit them.'' Important textbook treatments of
five-node count formulae discuss proof strategies but likewise omit
complete proofs. For example, \cite[p.~79]{estrada11} state that
``These calculations are based on the idea of spectral moments
analysed in the previous section,'' while
\cite[p.~136]{estrada_knight15} note that ``Formulae for a number of
simple subgraphs can be derived using very similar techniques to the
ones we have encountered so far.'' We hope that our main result and
its proof help to fill this gap.\\
\\
\noindent Both \cite{estrada11} and \cite{estrada_knight15} collect count
formulae for eight of the twenty-one connected subgraphs on five
nodes: in our notation, the cricket $|M_{79}^{(5)}|$, bull
$|M_{87}^{(5)}|$, banner $|M_{94}^{(5)}|$, lollipop
$|M_{117}^{(5)}|$, chevron $|M_{223}^{(5)}|$, hourglass
$|M_{235}^{(5)}|$, 5-circle $|M_{236}^{(5)}|$, and house
$|M_{237}^{(5)}|$. See
\cite[$|S_{8}|$ (4.18), $|S_{10}|$ (4.20), $|S_{9}|$ (4.19),
$|S_{11}|$ (4.21), $|S_{16}|$ (4.26), $|S_{12}|$ (4.22),
$|S_{7}|$ (4.17), $|S_{13}|$ (4.23)]{estrada11}, and
\cite[$|F_{9}|$, $|F_{11}|$, $|F_{10}|$, $|F_{12}|$ (p.~138),
$|F_{17}|$ (p.~139), $|F_{13}|$, $|F_{8}|$ (p.~138),
$|F_{14}|$ (p.~139)]{estrada_knight15}. The 5-star count
$|M_{75}^{(5)}|$ is also included in
\cite[$|S_{1,4}|$ (13.4)]{estrada_knight15}, and follows immediately
from the 4-star counts in \cite{alon_etal97,estrada11}.\\
\\
\noindent The same eight five-node subgraphs are treated in
\cite{alon_etal97}, apart from a known typographical omission in the
banner formula. In \cite[{$n_G(H_9)$}]{alon_etal97}, the term
$-2n_G(H_6)$ is missing; this term corrects twice for the unwanted
diamond subgraph. The authors note this error in the erratum available
from Noga Alon's webpage:
\url{https://web.math.princeton.edu/~nalon/PDFS/erratum1.pdf}.
There is also an unnecessary term in
\cite[$|S_{13}|$ (4.23)]{estrada11}, although the formula itself
remains correct: the term
$\mathbf{Q}\circ(\mathbf{A}^{3}\circ\mathbf{A})$, where
$\mathbf{Q}=\mathbf{A}^{2}\circ\mathbf{A}$, can be replaced by
$\mathbf{A}\circ\mathbf{A}^{2}\circ\mathbf{A}^{3}$, since
$\mathbf{A}\circ\mathbf{A}=\mathbf{A}$.\\
\\
\noindent The techniques and intended results in
\cite{estrada11,estrada_knight15} are closely related to those in
\cite{alon_etal97}, apart from differences in notation; in particular,
\cite{estrada11} favours a matrix formulation. Nevertheless, three of
the formulae presented in \cite{estrada11,estrada_knight15} are
incorrect because of (a) the omission of a correction term, (b) a
counting strategy that, unexpectedly, does not account for any
correction terms, or (c) incorrect indices of summation.\\
\\
\noindent \textbf{(a) Banner count $|S_{9}|$} in \cite[(4.19)]{estrada11} (as stated, this is correct only in a special case):
\begin{equation*}
    \text{``}
    |S_{9}| =
    \langle \mathbf{1} |
    \mathbf{P}-\mathbf{P}^{(D)}
    | \mathbf{k}'' \rangle;
    \quad
    \mathbf{P}
    =
    \frac{1}{2}
    [\mathbf{A}^{2}\circ(\mathbf{A}^{2}-1)];
    \quad
    \mathbf{k}''
    =
    \mathbf{k}-2;
    \quad
    \mathbf{k}
    =
    \mathbf{A}\mathbf{1}
    \text{''}.
\end{equation*}
Here $\langle \cdot|\cdot\rangle$ denotes the corresponding matrix
product, $\mathbf{1}$ is a column vector of ones, $\mathbf{P}^{(D)}$ is the diagonal matrix with diagonal entries $(\mathbf{P})_{ii}$, $\mathbf{A}$ is the adjacency matrix, $\mathbf{k}$ is the column vector of node degrees, and $\circ$ denotes the Hadamard product. Writing the scalar subtractions in dimensionally explicit form gives
\[
    \mathbf{P} = \frac{1}{2}
    \left[
        \mathbf{A}^{2}
        \circ
        (\mathbf{A}^{2}-\mathbf{1}\mathbf{1}')
    \right],
    \quad
    \mathbf{k}'' = \mathbf{k}-2\mathbf{1}.
\]
Thus $P_{ij} = \frac{1}{2}(g^2)_{ij}\bigl((g^2)_{ij}-1\bigr)$ and $ k_i'' = k_i-2.$ It follows that
\begin{align*}
    |S_{9}|
    &=
    \langle \mathbf{1}|\mathbf{P}|\mathbf{k}''\rangle
    -
    \langle \mathbf{1}|\mathbf{P}^{(D)}|\mathbf{k}''\rangle \\
    &=
    \frac{1}{2}
    \sum_{i,j}
    (g^2)_{ij}\bigl((g^2)_{ij}-1\bigr)(k_j-2) \\
    &\quad
    -
    \frac{1}{2}
    \sum_i
    (g^2)_{ii}\bigl((g^2)_{ii}-1\bigr)(k_i-2) \\
    &=
    \frac{1}{2}
    \sum_{i\neq j}
    (g^2)_{ij}\bigl((g^2)_{ij}-1\bigr)(k_j-2).
\end{align*}
Since $(g^2)_{ij}=(g^2)_{ji}$ and the sum is over all ordered pairs
$i\neq j$, we may equivalently write
\begin{equation*}
    |S_{9}| = \frac{1}{2}\sum_{\substack{i\neq j\\k_i>2}}
    (g^2)_{ij}\bigl((g^2)_{ij}-1\bigr)(k_i-2).
\end{equation*}
Comparing this with our banner count
\begin{equation*}
    |M_{94}^{(5)}| = \frac{1}{2}\sum_{\substack{i\neq j\\k_i>2}}
    (g^2)_{ij}\bigl((g^2)_{ij}-1\bigr)(k_i-2) - 2|M_{31}^{(4)}|,
\end{equation*}
it follows immediately that $|S_{9}|=|M_{94}^{(5)}|$ if and only if
$|M_{31}^{(4)}|=0$. Thus the stated formula is correct precisely when
$G$ contains no diamond subgraphs, and is not correct in general.\\
\\
\noindent \textbf{(b) Lollipop count $|S_{11}|$} in \cite[(4.21)]{estrada11} (as stated, this is correct only for triangle-free graphs):
\begin{equation*}
    \text{``} |S_{11}| = \frac{1}{2}
    \langle \mathbf{1} | \mathbf{A}^{2}-\mathbf{A}^{2(D)} | \mathbf{t} \rangle; \quad \mathbf{t} =
    \mathrm{diag}(\mathbf{A}^{3})\text{''},
\end{equation*}
where $\mathbf{1}$ is a column vector of ones, $\mathbf{A}$ is an
adjacency matrix, $\mathbf{A}^{2(D)}$ is the diagonal matrix with
diagonal entries $(\mathbf{A}^{2})_{ii}$, and $\mathbf{t}$ is the
vector with entries $(\mathbf{A}^{3})_{ii}$. We can write
\begin{equation*}
    |S_{11}| = \frac{1}{2}
    \langle \mathbf{1}|\mathbf{A}^{2}|\mathbf{t}\rangle -
    \frac{1}{2} \langle \mathbf{1}|\mathbf{A}^{2(D)}|\mathbf{t}\rangle.
\end{equation*}
The apparent strategy underlying the formula in \cite{estrada11} is
to consider paths of length 2 between distinct nodes $i$ and $j$,
where one endpoint is also incident to a triangle, and to remove the
diagonal case $i=j$. However, this does not correct for the other
unwanted configurations that arise from overlaps between the
length-2 path and the triangle. In our notation, and using the symmetry of $g^2$ to interchange the indices if necessary, the expression becomes
\begin{align*}
    |S_{11}|
    &=
    \frac{1}{2}
    \sum_{i,j}(g^2)_{ij}(g^3)_{ii}
    -
    \frac{1}{2}
    \sum_i(g^2)_{ii}(g^3)_{ii} \\
    &=
    \frac{1}{2}
    \sum_{i\neq j}(g^2)_{ij}(g^3)_{ii}.
\end{align*}
For completeness,
\begin{align*}
    \frac{1}{2}
    \sum_i(g^2)_{ii}(g^3)_{ii}
    &=
    \frac{1}{2}
    \sum_i k_i(g^3)_{ii} \\
    &=
    \frac{1}{2}
    \sum_i(g^3)_{ii}(k_i-2)
    +
    \tr(g^3) \\
    &=
    |M_{15}^{(4)}|
    +
    6|M_{7}^{(3)}|,
\end{align*}
using \cite[eqns.~(6) and~(9)]{agasse-duval_lawford25}, although this
identity is not needed because the diagonal terms cancel above. Expanding $(g^2)_{ij}$ gives
\begin{align*}
    |S_{11}|
    &=
    \frac{1}{2}\sum_{i\neq j}\sum_r g_{ir}g_{rj}(g^3)_{ii} \\
    &=
    \frac{1}{2} \sum_{\substack{(i,r):\,\{i,r\}\in E}}
    (g^3)_{ii}(k_r-1),
\end{align*}
where the final sum is over both orientations of each edge. Comparing
this with our lollipop count,
\begin{equation*}
    |M_{117}^{(5)}| = \frac{1}{2}\sum_{\substack{(i,r):\,\{i,r\}\in E}}(g^3)_{ii}(k_r-1) - 6|M_{7}^{(3)}| - 2|M_{15}^{(4)}|
    - 4|M_{31}^{(4)}|,
\end{equation*}
we obtain
\begin{equation*}
    |S_{11}|-|M_{117}^{(5)}| = 6|M_{7}^{(3)}| +
    2|M_{15}^{(4)}| + 4|M_{31}^{(4)}|.
\end{equation*}
Since all subgraph counts are non-negative, $|S_{11}|=|M_{117}^{(5)}|$ if and only if
\[
    |M_{7}^{(3)}|=0.
\]
Thus the stated formula is correct precisely for triangle-free graphs,
in which case both sides are necessarily zero, and is not correct in
general.\\
\\
\noindent \textbf{(c) Banner count $|F_{10}|$} in
\cite[p.~138]{estrada_knight15} (as stated, this is not correct in
general):
\begin{equation*}
    \text{``} |F_{10}| =
    \frac{1}{2}\sum_i(k_i-2) \times
    \sum_{i,j}\binom{(A^2)_{ij}}{2} - 2|F_7|
    \text{''},
\end{equation*}
where $k_i$ is the degree of node $i$, $A$ is an adjacency matrix,
and $|F_7|$ is the diamond count. Interpreted literally, the two sums
in the product are independent. By contrast, our banner count is
\begin{equation*}
    |M_{94}^{(5)}| =
    \sum_{\substack{i\neq j\\k_i>2}}
    \binom{(g^2)_{ij}}{2}(k_i-2) - 2|M_{31}^{(4)}|,
\end{equation*}
or, equivalently,
\begin{equation*}
    |M_{94}^{(5)}| = \sum_{\substack{i\\k_i>2}}(k_i-2)
    \sum_{j\neq i}\binom{(g^2)_{ij}}{2} - 2|M_{31}^{(4)}|.
\end{equation*}
Thus the degree factor $(k_i-2)$ must be associated with the same
node $i$ that indexes the corresponding row of $g^2$; it cannot be
factored out into an independent sum over all nodes. A direct counterexample is given by taking $G$ itself to be a banner.
Its degree sequence is $(3,2,2,2,1)$, so that
\[
    \sum_i(k_i-2)=0,
\]
and it contains no diamond subgraph. Hence the formula stated in
\cite{estrada_knight15} gives $|F_{10}|=0$. However, by construction,
\[
    |M_{94}^{(5)}|=1.
\]
The discrepancy therefore appears to be an indexing error. A
corrected version, in the notation of \cite{estrada_knight15}, is
\begin{equation*}
    |F_{10}| = \sum_{\substack{i\\k_i>2}}(k_i-2)
    \sum_{j\neq i} \binom{(A^2)_{ij}}{2} - 2|F_7|.
\end{equation*}
    
\end{example}

\begin{example}[An empirical application to real-world network data]
\label{example:counts_WN_2013_4}

We now demonstrate that the analytic formulae of
Theorem~\ref{thm:analytic_count_5node} can be implemented directly in
an empirical application. The implementations are intended to
illustrate and verify the formulae rather than to provide optimized
algorithms for large-scale subgraph enumeration, and we therefore do
not undertake runtime comparisons with specialized exact or sampling
methods. There is a substantial literature in applied mathematics and computer science on efficient exact and approximate subgraph enumeration. For five-node subgraph counting, \cite{pinar_etal17} develop algorithms
based on decomposing subgraphs into a small set of smaller patterns
that can be efficiently enumerated. Related work, including methods
applicable to large networks, includes
\cite{bera_etal20,chen_lui18,curticapean_etal17,hocevar_demsar16,
kane_etal12,pashanasangi_seshadhri20}; see also the survey by
\cite{ribeiro_etal21}. Our objective here is different: to illustrate
the direct implementation and empirical use of the analytical
formulae derived above.\\
\\
\noindent Using Python 3, we implemented the count formula for each five-node subgraph appearing in Table~\ref{tab:subgraphs5}. We applied these
routines to a small real-world airline network with 88 nodes and 522
edges, described in detail by Agasse-Duval and Lawford
\cite{agasse-duval_lawford25}, who examine three-node and four-node
motifs in airline networks, and by Lawford and Mehmeti
\cite{lawford_mehmeti20}, who use analytic five-node subgraph counts
for the 5-star, 5-path, 5-arrow, kite, and 5-complete subgraphs to
compute generalized clustering coefficients of order five. We then
computed each corresponding induced count using
Theorem~\ref{thm:analytic_count_not_nested_5node}.\\
\\
\noindent To provide a more informative comparison of the induced five-node subgraph profile, we use two null ensembles, each containing 1,000
networks. The first consists of Erd{\H{o}}s--R{\'e}nyi graphs
$G(n,p)$ with $n=88$ and $p=522/\binom{88}{2}\simeq 0.13636$, so that the number of nodes and expected edge density agree with the observed network. The second consists of degree-preserving randomizations of the observed network, generated by repeated double-edge swaps. Each realization starts from the observed graph and undergoes $50|E|=26{,}100$ successful swaps, while self-loops and multiple edges are excluded. Thus the complete observed degree sequence is preserved exactly. For each induced five-node subgraph $M_a^{(5)}$, let $\mu_a$ and $s_a$ denote the mean and sample standard deviation of its induced count
under a given null ensemble. We summarize the comparison by the
standardized deviation $Z_a = (|\widetilde M_a^{(5)}|-\mu_a)/s_a$. Figure~\ref{fig:WN_2013_4_subgraph_counts_main}(a) displays the general
and induced counts in the observed Southwest Airlines network, while
panels~(b) and~(c) give the standardized induced-count profiles relative to the Erd{\H{o}}s--R{\'e}nyi and degree-preserving null ensembles.\\
\\
\noindent The distinction between the two null models is substantial. Relative to the Erd{\H{o}}s--R{\'e}nyi ensemble, many five-node subgraphs show very large standardized deviations. Much of this apparent structure,
however, disappears once the observed degree sequence is preserved. For example, the 5-complete has $Z^{\mathrm{ER}}\simeq 11{,}514$ but
$Z^{\mathrm{deg}}\simeq 0.34$, while the corresponding values for the
lamp are approximately $1{,}093$ and $-0.08$, and for the 5-star
approximately $176$ and $0.77$. These comparisons indicate that much
of the departure from an Erd{\H{o}}s--R{\'e}nyi graph is attributable
to the heterogeneous degree sequence of the airline network.\\
\\
\noindent Some deviations remain after controlling for the degree sequence. In particular, the crown is under-represented
($Z^{\mathrm{deg}}\simeq -6.04$), while the hourglass and chevron are
over-represented ($Z^{\mathrm{deg}}\simeq 4.56$ and $3.96$, respectively); the cricket and 5-circle are also under-represented
($Z^{\mathrm{deg}}\simeq -3.55$ and $-3.14$). We interpret these
values descriptively as standardized differences from the null
ensembles rather than as a formal multiple-testing analysis of
network motifs. Nevertheless, they illustrate how induced five-node
subgraph profiles can distinguish higher-order structural features
from effects that are largely explained by the degree sequence. Such differences may be useful for characterizing real-world networks
and for applications such as link prediction, as also noted by
\cite{pinar_etal17}. In Table~\ref{tab:WN_2013_4_subgraph_examples}, we show representative subgraphs together with the corresponding node identities. Airport identities associated with the three-letter IATA codes can be checked using the IATA code-search directory:
\url{https://www.iata.org/en/publications/directories/code-search/}.

\begin{figure}[p]
    \centering
    \includegraphics[width=\linewidth,height=0.72\textheight,keepaspectratio]{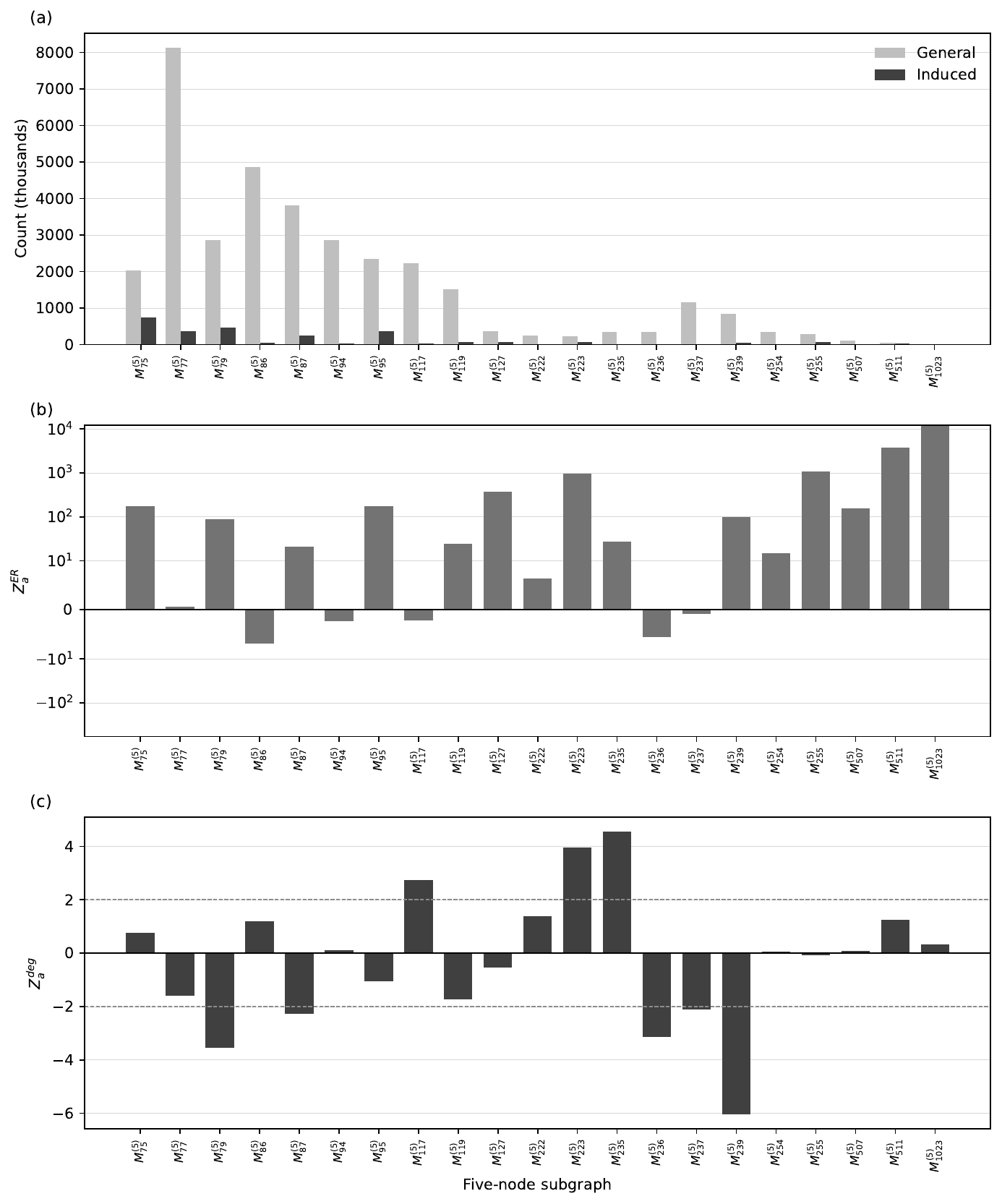}
    \caption{Empirical five-node subgraph counts and null-model comparisons for the Southwest Airlines route network in 2013Q4.
    Panel~(a) shows the general five-node subgraph counts from
    Theorem~\ref{thm:analytic_count_5node} and the corresponding induced
    counts obtained from Theorem~\ref{thm:analytic_count_not_nested_5node}.
    Panel~(b) shows standardized deviations of the induced counts relative
    to an ensemble of 1,000 Erd{\H{o}}s--R{\'e}nyi graphs $G(n,p)$ with
    $n=88$ and $p=522/\binom{88}{2}\simeq0.13636$; the vertical axis uses
    a symmetric logarithmic scale.
    Panel~(c) shows the corresponding standardized deviations relative to
    1,000 degree-preserving randomized networks generated by double-edge
    swaps. The dashed horizontal lines in panel~(c) indicate $Z=\pm2$.
    See Example~\ref{example:counts_WN_2013_4} for details.}
    \label{fig:WN_2013_4_subgraph_counts_main}
\end{figure}

\end{example}

\clearpage

\begin{landscape}
\begin{table}[p]
\centering
\includegraphics[width=0.96\linewidth,height=0.82\textheight,keepaspectratio]{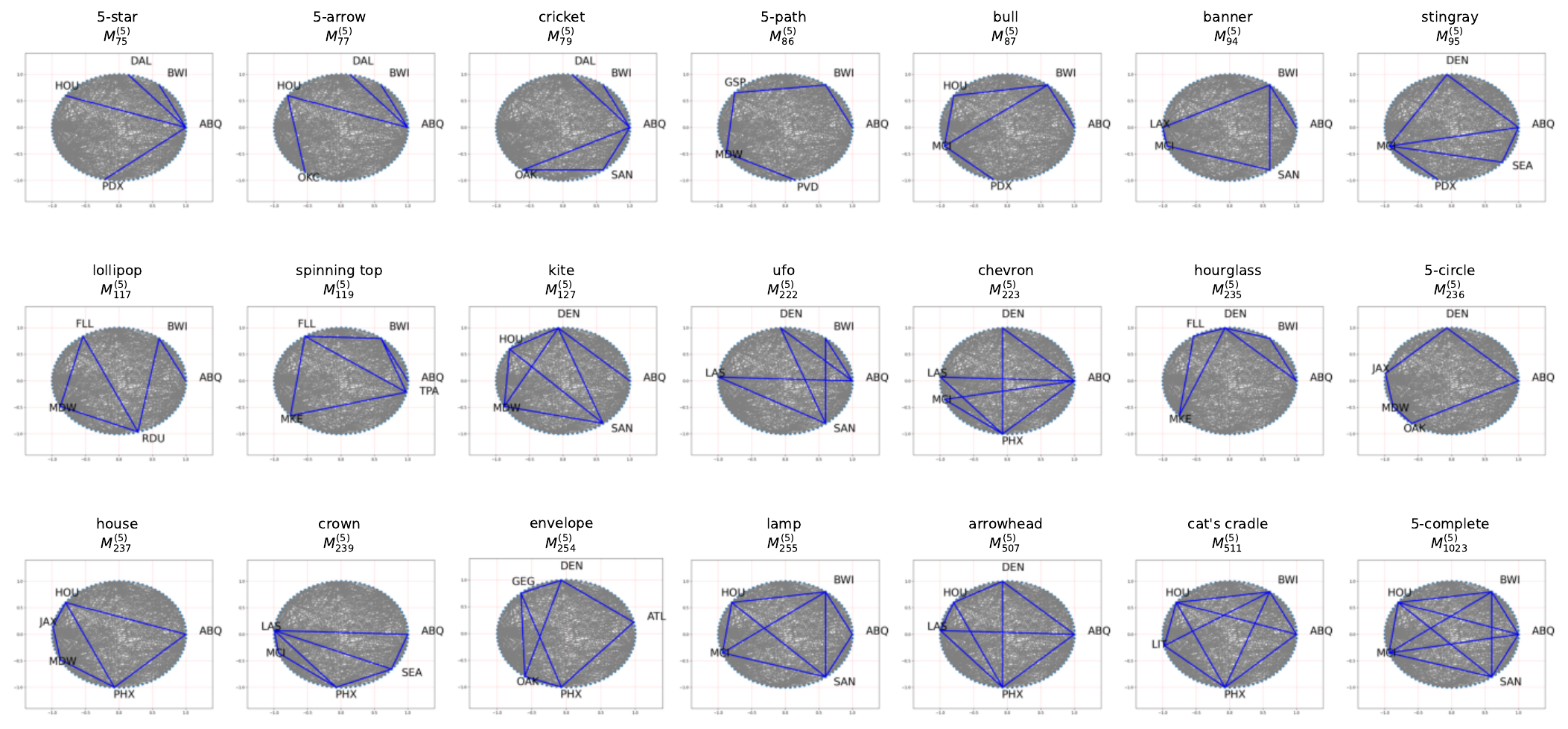}
\caption{Representative five-node subgraphs from the Southwest Airlines route network in 2013Q4, as discussed in Example~\ref{example:counts_WN_2013_4}. The edges of each selected
five-node subgraph are highlighted, with the remainder of the network
shown in the background. Airport identities corresponding to the
three-letter IATA codes can be found in the IATA code-search directory.}
\label{tab:WN_2013_4_subgraph_examples}
\end{table}
\end{landscape}

\section{Conclusions}\label{sec:conclusions}
Exact subgraph count formulae can lead to a better understanding of the topological structure and related function of real-world networks, and can provide improved insight into widely used graph statistics such as the assortativity index and the overall clustering coefficient. Much of the focus of research in applied mathematics since the 1970s and, more recently, in theoretical and applied computer science has been on the development of faster implementations of such formulae and their application to ever larger datasets. This is of deep theoretical interest and crucial empirical importance. However, there has been relatively little work on elementary combinatorial treatments of small subgraph counts. While some of these formulae and proof strategies now appear in important textbooks on graph theory, their correct derivation can be surprisingly tricky, and new formulations of exact count formulae can be conveniently specialized to derive non-trivial properties of various theoretical graphs.\\
\\
\noindent Here we have presented general and induced counts for all 21 connected subgraphs on five nodes, several of whose general-count formulae appear to be new. We place particular emphasis on the combinatorial intuition behind the results and illustrate, through a series of examples, how these formulae can be used in theoretical and empirical work. In the airline application, degree-preserving randomization shows that many of the large departures from an Erd{\H{o}}s--R{\'e}nyi benchmark are explained by degree heterogeneity, while several higher-order deviations remain. While exact formulae might, in future, motivate useful computational improvements, we expect they will be most useful in deriving theoretical results for special graphs and graph statistics.\\
\\
\noindent The neighbourhood-subgraph construction also admits a general lifting identity: whenever a subgraph pattern has a distinguished node adjacent to every other node, its count can be obtained from counts of the remaining pattern inside node neighbourhoods, with a correction for structurally equivalent choices of the distinguished node. This identity underlies several of the five-node formulae and extends immediately to suitable patterns of larger order. More generally, the constructions developed here can yield exact counts for individual larger subgraphs, as illustrated by the six-node examples in Remark~\ref{rem:results_strategies_extensions}\ref{rem:six_node_examples}. They do not, however, provide a closed recursive scheme for complete
enumeration at larger orders, since the required correction terms need
not themselves involve only simpler patterns. A further direction is
to combine exact subgraph count formulae more directly with
probabilistic network models; for example, analytic expressions for
the expectations and higher moments of five-node counts under standard
random-graph models could provide useful benchmarks for formal
comparisons of observed network structure. It would also be
interesting to investigate whether parts of the derivation process for
larger subgraphs, in particular, the identification of candidate
constructions and their correction terms, can be systematized or
partially automated while retaining a transparent combinatorial
interpretation. Developing more systematic methods that overcome this
limitation remains an interesting direction for future work.

\section*{Funding}

This research received no specific grant from any funding agency in
the public, commercial, or not-for-profit sectors.

\section*{Acknowledgements}

The author is grateful to Jack Lawford and Ian Wanless for helpful
comments and discussions.

\appendix

\counterwithin{table}{section}
\counterwithin{figure}{section}

\newpage

\section{Counting induced subgraphs}\label{sec:induced_counts}

Induced subgraph counts $\widetilde{y}=(|\widetilde{M}_{a}^{(b)}|)$ can be recovered from the corresponding general subgraph counts $y=(|M_{a}^{(b)}|)$ through the linear relation $\widetilde{y}=y-A\widetilde{y}$. Equivalently, $\widetilde{y}=(I+A)^{-1}y$. Under the ordering of subgraphs used here, $A$ is strictly upper triangular, so that $I+A$ is unit upper triangular and hence invertible. The matrix $A$ for five-node subgraphs is given in Table~\ref{tab:matrix_A}. We give the resulting individual induced count formulae in
Theorem~\ref{thm:analytic_count_not_nested_5node}.

\begin{theorem}[Count formulae for induced subgraphs on five nodes]
\label{thm:analytic_count_not_nested_5node}
\begin{align*}
|\widetilde{M}_{75}^{(5)}|
={}& |M_{75}^{(5)}| - |M_{79}^{(5)}| + |M_{95}^{(5)}|
     - |M_{127}^{(5)}| - 2|M_{223}^{(5)}| \\
&+ |M_{235}^{(5)}| - |M_{239}^{(5)}| + 2|M_{255}^{(5)}|
  + |M_{507}^{(5)}| - 3|M_{511}^{(5)}|
  + 5|M_{1023}^{(5)}|. \\[0.5ex]
|\widetilde{M}_{77}^{(5)}|
={}& |M_{77}^{(5)}| - 2|M_{79}^{(5)}| - 2|M_{87}^{(5)}|
     - 2|M_{94}^{(5)}| + 5|M_{95}^{(5)}| - |M_{117}^{(5)}| \\
&+ 4|M_{119}^{(5)}| - 9|M_{127}^{(5)}| + 6|M_{222}^{(5)}|
  - 12|M_{223}^{(5)}| + 4|M_{235}^{(5)}| + 4|M_{237}^{(5)}| \\
&- 10|M_{239}^{(5)}| - 10|M_{254}^{(5)}| + 20|M_{255}^{(5)}|
  + 20|M_{507}^{(5)}| - 36|M_{511}^{(5)}|
  + 60|M_{1023}^{(5)}|. \\[0.5ex]
|\widetilde{M}_{79}^{(5)}|
={}& |M_{79}^{(5)}| - 2|M_{95}^{(5)}| + 3|M_{127}^{(5)}|
     + 6|M_{223}^{(5)}| - 2|M_{235}^{(5)}| \\
&+ 3|M_{239}^{(5)}| - 8|M_{255}^{(5)}| - 4|M_{507}^{(5)}|
  + 15|M_{511}^{(5)}| - 30|M_{1023}^{(5)}|. \\[0.5ex]
|\widetilde{M}_{86}^{(5)}|
={}& |M_{86}^{(5)}| - |M_{87}^{(5)}| - 2|M_{94}^{(5)}|
     + 2|M_{95}^{(5)}| - 2|M_{117}^{(5)}| + 4|M_{119}^{(5)}| \\
&- 6|M_{127}^{(5)}| + 6|M_{222}^{(5)}| - 6|M_{223}^{(5)}|
  + 4|M_{235}^{(5)}| - 5|M_{236}^{(5)}| + 7|M_{237}^{(5)}| \\
&- 10|M_{239}^{(5)}| - 14|M_{254}^{(5)}| + 18|M_{255}^{(5)}|
  + 24|M_{507}^{(5)}| - 36|M_{511}^{(5)}|
  + 60|M_{1023}^{(5)}|. \\[0.5ex]
|\widetilde{M}_{87}^{(5)}|
={}& |M_{87}^{(5)}| - 2|M_{95}^{(5)}| - 2|M_{119}^{(5)}|
     + 6|M_{127}^{(5)}| + 6|M_{223}^{(5)}| - |M_{237}^{(5)}| \\
&+ 5|M_{239}^{(5)}| + 4|M_{254}^{(5)}| - 14|M_{255}^{(5)}|
  - 12|M_{507}^{(5)}| + 30|M_{511}^{(5)}|
  - 60|M_{1023}^{(5)}|. \\[0.5ex]
|\widetilde{M}_{94}^{(5)}|
={}& |M_{94}^{(5)}| - |M_{95}^{(5)}| - |M_{119}^{(5)}|
     + 3|M_{127}^{(5)}| - 6|M_{222}^{(5)}| + 6|M_{223}^{(5)}| \\
&- 2|M_{237}^{(5)}| + 4|M_{239}^{(5)}| + 8|M_{254}^{(5)}|
  - 12|M_{255}^{(5)}| - 16|M_{507}^{(5)}|
  + 30|M_{511}^{(5)}| - 60|M_{1023}^{(5)}|. \\[0.5ex]
|\widetilde{M}_{95}^{(5)}|
={}& |M_{95}^{(5)}| - 3|M_{127}^{(5)}| - 6|M_{223}^{(5)}|
     - 2|M_{239}^{(5)}| + 10|M_{255}^{(5)}| \\
&+ 4|M_{507}^{(5)}| - 24|M_{511}^{(5)}|
  + 60|M_{1023}^{(5)}|. \\[0.5ex]
|\widetilde{M}_{117}^{(5)}|
={}& |M_{117}^{(5)}| - 2|M_{119}^{(5)}| + 3|M_{127}^{(5)}|
     - 4|M_{235}^{(5)}| - 2|M_{237}^{(5)}| + 6|M_{239}^{(5)}| \\
&+ 6|M_{254}^{(5)}| - 12|M_{255}^{(5)}| - 16|M_{507}^{(5)}|
  + 30|M_{511}^{(5)}| - 60|M_{1023}^{(5)}|. \\[0.5ex]
|\widetilde{M}_{119}^{(5)}|
={}& |M_{119}^{(5)}| - 3|M_{127}^{(5)}| - 2|M_{239}^{(5)}|
     - 2|M_{254}^{(5)}| + 8|M_{255}^{(5)}| \\
&+ 8|M_{507}^{(5)}| - 24|M_{511}^{(5)}|
  + 60|M_{1023}^{(5)}|. \\[0.5ex]
|\widetilde{M}_{127}^{(5)}|
={}& |M_{127}^{(5)}| - 2|M_{255}^{(5)}|
     + 6|M_{511}^{(5)}| - 20|M_{1023}^{(5)}|. \\[0.5ex]
|\widetilde{M}_{222}^{(5)}|
={}& |M_{222}^{(5)}| - |M_{223}^{(5)}| - |M_{254}^{(5)}|
     + |M_{255}^{(5)}| + 2|M_{507}^{(5)}| \\
&- 4|M_{511}^{(5)}| + 10|M_{1023}^{(5)}|. \\[0.5ex]
|\widetilde{M}_{223}^{(5)}|
={}& |M_{223}^{(5)}| - |M_{255}^{(5)}|
     + 3|M_{511}^{(5)}| - 10|M_{1023}^{(5)}|. \\[0.5ex]
|\widetilde{M}_{235}^{(5)}|
={}& |M_{235}^{(5)}| - |M_{239}^{(5)}| + 2|M_{255}^{(5)}|
     + 2|M_{507}^{(5)}| - 6|M_{511}^{(5)}|
     + 15|M_{1023}^{(5)}|. \\[0.5ex]
|\widetilde{M}_{236}^{(5)}|
={}& |M_{236}^{(5)}| - |M_{237}^{(5)}| + |M_{239}^{(5)}|
     + 2|M_{254}^{(5)}| - 2|M_{255}^{(5)}| \\
&- 4|M_{507}^{(5)}| + 6|M_{511}^{(5)}|
  - 12|M_{1023}^{(5)}|.
\end{align*}
\pagebreak
\begin{align*}
|\widetilde{M}_{237}^{(5)}|
={}& |M_{237}^{(5)}| - 2|M_{239}^{(5)}| - 4|M_{254}^{(5)}|
     + 6|M_{255}^{(5)}| + 12|M_{507}^{(5)}| \\
&- 24|M_{511}^{(5)}| + 60|M_{1023}^{(5)}|. \\[0.5ex]
|\widetilde{M}_{239}^{(5)}|
={}& |M_{239}^{(5)}| - 4|M_{255}^{(5)}| - 4|M_{507}^{(5)}|
     + 18|M_{511}^{(5)}| - 60|M_{1023}^{(5)}|. \\[0.5ex]
|\widetilde{M}_{254}^{(5)}|
={}& |M_{254}^{(5)}| - |M_{255}^{(5)}| - 4|M_{507}^{(5)}|
     + 9|M_{511}^{(5)}| - 30|M_{1023}^{(5)}|. \\[0.5ex]
|\widetilde{M}_{255}^{(5)}|
={}& |M_{255}^{(5)}| - 6|M_{511}^{(5)}|
     + 30|M_{1023}^{(5)}|. \\[0.5ex]
|\widetilde{M}_{507}^{(5)}|
={}& |M_{507}^{(5)}| - 3|M_{511}^{(5)}|
     + 15|M_{1023}^{(5)}|. \\[0.5ex]
|\widetilde{M}_{511}^{(5)}|
={}& |M_{511}^{(5)}| - 10|M_{1023}^{(5)}|. \\[0.5ex]
|\widetilde{M}_{1023}^{(5)}|
={}& |M_{1023}^{(5)}|.
\end{align*}
\end{theorem}

\begin{proof}
All of the results follow directly from $\widetilde{y}=(I+A)^{-1}y$.
\end{proof}

\noindent
To develop some combinatorial intuition, we consider the envelope count $|\widetilde{M}_{254}^{(5)}|$. In what follows, references to the degrees of nodes refer to their degrees within the corresponding
five-node pattern. There is one way to remove an edge from a lamp to
obtain an envelope, namely the edge between its two degree-4 nodes, and four ways to remove an edge from an arrowhead to obtain an envelope, each involving a degree-4 and a degree-3 node. Further, there are six ways to remove an edge from a cat's cradle to obtain a lamp, by removing an edge between a degree-3 and a degree-4 node, and three ways to remove an edge from a cat's cradle to obtain an arrowhead, by removing an edge between two degree-4 nodes. Hence there are $(1/2)((6\times1)+(3\times4))=9$ ways to remove two edges from a cat's cradle to obtain an envelope, where division by two corrects for the two possible orders in which the two edges can be removed. Finally, there are 30 ways to remove three edges from a 5-complete subgraph to obtain an envelope. To see this, choose one of the 10 edges of the 5-complete to be the isolated edge among the three removed edges. The remaining three nodes form a triangle, and removing any two of its three edges gives the other two removed edges. Thus there are $10\times3=30$ possibilities. Putting these observations together gives
\begin{equation*}
    |\widetilde{M}_{254}^{(5)}|
    =
    |M_{254}^{(5)}|
    -|\widetilde{M}_{255}^{(5)}|
    -4|\widetilde{M}_{507}^{(5)}|
    -9|\widetilde{M}_{511}^{(5)}|
    -30|\widetilde{M}_{1023}^{(5)}|.
\end{equation*}
Since
$|\widetilde{M}_{1023}^{(5)}|=|M_{1023}^{(5)}|$, substitution of the
preceding induced-count identities gives
\begin{equation*}
    |\widetilde{M}_{254}^{(5)}|
    =
    |M_{254}^{(5)}|
    -|M_{255}^{(5)}|
    -4|M_{507}^{(5)}|
    +9|M_{511}^{(5)}|
    -30|M_{1023}^{(5)}|,
\end{equation*}
as required. A full proof of
Theorem~\ref{thm:analytic_count_not_nested_5node} using this
combinatorial approach would be very cumbersome.

\begin{table}[p]
\centering
{\scriptsize
\renewcommand{\arraystretch}{1.65}
\setlength{\arraycolsep}{2.2pt}
\(
\begin{blockarray}{rccccccccccccccccccccc}
\begin{block}{r(ccccccccccccccccccccc)}
  \textrm{5-star} \; |M_{75}^{(5)}|
  & 0 & 0 & \hlight{1} & 0 & 0 & 0 & \hlight{1} & 0 & 0
  & \hlight{1} & 0 & \hlight{2} & \hlight{1} & 0 & 0 & \hlight{1}
  & 0 & \hlight{2} & \hlight{1} & \hlight{3} & \hlight{5}\\
  \textrm{5-arrow} \; |M_{77}^{(5)}|
  & 0 & 0 & \hlight{2} & 0 & \hlight{2} & \hlight{2} & \hlight{5}
  & \hlight{1} & \hlight{4} & \hlight{9} & \hlight{6} & \hlight{12}
  & \hlight{4} & 0 & \hlight{4} & \hlight{10} & \hlight{10}
  & \hlight{20} & \hlight{20} & \hlight{36} & \hlight{60}\\
  \textrm{cricket} \; |M_{79}^{(5)}|
  & 0 & 0 & 0 & 0 & 0 & 0 & \hlight{2} & 0 & 0 & \hlight{3}
  & 0 & \hlight{6} & \hlight{2} & 0 & 0 & \hlight{3} & 0
  & \hlight{8} & \hlight{4} & \hlight{15} & \hlight{30}\\
  \textrm{5-path} \; |M_{86}^{(5)}|
  & 0 & 0 & 0 & 0 & \hlight{1} & \hlight{2} & \hlight{2}
  & \hlight{2} & \hlight{4} & \hlight{6} & \hlight{6} & \hlight{6}
  & \hlight{4} & \hlight{5} & \hlight{7} & \hlight{10} & \hlight{14}
  & \hlight{18} & \hlight{24} & \hlight{36} & \hlight{60}\\
  \textrm{bull} \; |M_{87}^{(5)}|
  & 0 & 0 & 0 & 0 & 0 & 0 & \hlight{2} & 0 & \hlight{2}
  & \hlight{6} & 0 & \hlight{6} & 0 & 0 & \hlight{1}
  & \hlight{5} & \hlight{4} & \hlight{14} & \hlight{12}
  & \hlight{30} & \hlight{60}\\
  \textrm{banner} \; |M_{94}^{(5)}|
  & 0 & 0 & 0 & 0 & 0 & 0 & \hlight{1} & 0 & \hlight{1}
  & \hlight{3} & \hlight{6} & \hlight{6} & 0 & 0 & \hlight{2}
  & \hlight{4} & \hlight{8} & \hlight{12} & \hlight{16}
  & \hlight{30} & \hlight{60}\\
  \textrm{stingray} \; |M_{95}^{(5)}|
  & 0 & 0 & 0 & 0 & 0 & 0 & 0 & 0 & 0 & \hlight{3} & 0
  & \hlight{6} & 0 & 0 & 0 & \hlight{2} & 0 & \hlight{10}
  & \hlight{4} & \hlight{24} & \hlight{60}\\
  \textrm{lollipop} \; |M_{117}^{(5)}|
  & 0 & 0 & 0 & 0 & 0 & 0 & 0 & 0 & \hlight{2} & \hlight{3}
  & 0 & 0 & \hlight{4} & 0 & \hlight{2} & \hlight{6}
  & \hlight{6} & \hlight{12} & \hlight{16} & \hlight{30}
  & \hlight{60}\\
  \textrm{spinning top} \; |M_{119}^{(5)}|
  & 0 & 0 & 0 & 0 & 0 & 0 & 0 & 0 & 0 & \hlight{3}
  & 0 & 0 & 0 & 0 & 0 & \hlight{2} & \hlight{2}
  & \hlight{8} & \hlight{8} & \hlight{24} & \hlight{60}\\
  \textrm{kite} \; |M_{127}^{(5)}|
  & 0 & 0 & 0 & 0 & 0 & 0 & 0 & 0 & 0 & 0 & 0 & 0 & 0 & 0
  & 0 & 0 & 0 & \hlight{2} & 0 & \hlight{6} & \hlight{20}\\
  \textrm{ufo} \; |M_{222}^{(5)}|
  & 0 & 0 & 0 & 0 & 0 & 0 & 0 & 0 & 0 & 0 & 0 & \hlight{1}
  & 0 & 0 & 0 & 0 & \hlight{1} & \hlight{1} & \hlight{2}
  & \hlight{4} & \hlight{10}\\
  \textrm{chevron} \; |M_{223}^{(5)}|
  & 0 & 0 & 0 & 0 & 0 & 0 & 0 & 0 & 0 & 0 & 0 & 0 & 0 & 0
  & 0 & 0 & 0 & \hlight{1} & 0 & \hlight{3} & \hlight{10}\\
  \textrm{hourglass} \; |M_{235}^{(5)}|
  & 0 & 0 & 0 & 0 & 0 & 0 & 0 & 0 & 0 & 0 & 0 & 0 & 0 & 0
  & 0 & \hlight{1} & 0 & \hlight{2} & \hlight{2}
  & \hlight{6} & \hlight{15}\\
  \textrm{5-circle} \; |M_{236}^{(5)}|
  & 0 & 0 & 0 & 0 & 0 & 0 & 0 & 0 & 0 & 0 & 0 & 0 & 0 & 0
  & \hlight{1} & \hlight{1} & \hlight{2} & \hlight{2}
  & \hlight{4} & \hlight{6} & \hlight{12}\\
  \textrm{house} \; |M_{237}^{(5)}|
  & 0 & 0 & 0 & 0 & 0 & 0 & 0 & 0 & 0 & 0 & 0 & 0 & 0 & 0 & 0
  & \hlight{2} & \hlight{4} & \hlight{6} & \hlight{12}
  & \hlight{24} & \hlight{60}\\
  \textrm{crown} \; |M_{239}^{(5)}|
  & 0 & 0 & 0 & 0 & 0 & 0 & 0 & 0 & 0 & 0 & 0 & 0 & 0 & 0 & 0
  & 0 & 0 & \hlight{4} & \hlight{4} & \hlight{18} & \hlight{60}\\
  \textrm{envelope} \; |M_{254}^{(5)}|
  & 0 & 0 & 0 & 0 & 0 & 0 & 0 & 0 & 0 & 0 & 0 & 0 & 0 & 0 & 0
  & 0 & 0 & \hlight{1} & \hlight{4} & \hlight{9} & \hlight{30}\\
  \textrm{lamp} \; |M_{255}^{(5)}|
  & 0 & 0 & 0 & 0 & 0 & 0 & 0 & 0 & 0 & 0 & 0 & 0 & 0 & 0 & 0
  & 0 & 0 & 0 & 0 & \hlight{6} & \hlight{30}\\
  \textrm{arrowhead} \; |M_{507}^{(5)}|
  & 0 & 0 & 0 & 0 & 0 & 0 & 0 & 0 & 0 & 0 & 0 & 0 & 0 & 0 & 0
  & 0 & 0 & 0 & 0 & \hlight{3} & \hlight{15}\\
  \textrm{cat's cradle} \; |M_{511}^{(5)}|
  & 0 & 0 & 0 & 0 & 0 & 0 & 0 & 0 & 0 & 0 & 0 & 0 & 0 & 0 & 0
  & 0 & 0 & 0 & 0 & 0 & \hlight{10}\\
  \textrm{5-complete} \; |M_{1023}^{(5)}|
  & 0 & 0 & 0 & 0 & 0 & 0 & 0 & 0 & 0 & 0 & 0 & 0 & 0 & 0 & 0
  & 0 & 0 & 0 & 0 & 0 & 0\\
\end{block}
\end{blockarray}
\)
}
\caption{Matrix $A$, whose off-diagonal entry in row $r$ and column $s$ is the number of copies of the subgraph corresponding to row $r$
contained in the subgraph corresponding to column $s$; rows and columns use the same ordering. Up to a simultaneous reordering of rows and columns, $I+A$ appears in Figure~12 of the working-paper version of \cite{pinar_etal17}. See
Theorem~\ref{thm:analytic_count_not_nested_5node} and
Remark~\ref{remark:matrix_A}.}
\label{tab:matrix_A}
\end{table}

\begin{remark}\label{remark:matrix_A}
With regard to the construction of the matrix $A$ in Table~\ref{tab:matrix_A}, we note that:
\begin{enumerate}[label=(\roman*)]
    \item \label{1}
    The induced subgraph counts are of the form
    \[
        |\widetilde{M}_{a}^{(b)}| = |M_{a}^{(b)}| -
        \sum_{s>a} c_{a,s}\,|\widetilde{M}_{s}^{(b)}|,
    \]
    where the coefficients $c_{a,s}$ are entries of $A$; that is,
    \[
        \textrm{induced count} =
        \textrm{general subgraph count} - \textrm{correction}.
    \]
    No induced subgraph $\widetilde{M}_{s}^{(b)}$ with $s\leq a$ can
    appear in the correction term. To see this, recall that the
    subgraphs are ordered by increasing $a$, where $a$ is the smallest
    binary representation of the subgraph's adjacency matrix over all
    node labellings (see Section~\ref{sec:notation}). Suppose that
    $M_{a}^{(b)}$ is obtained from $M_{s}^{(b)}$ by deleting at least
    one edge. Choose a labelling of $M_{s}^{(b)}$ that attains its
    canonical code $s$. Under the same labelling, deleting an edge
    changes a binary digit from 1 to 0 and therefore strictly decreases the corresponding binary value. Since $a$ is the minimum over all labellings of the resulting graph, it follows that $a<s$. Consequently, the main diagonal and strictly lower triangular part of $A$ are zero.\\
    \\
    \noindent For example, deleting one edge from a stingray
    $\widetilde{M}_{95}^{(5)}$, which has six edges, cannot produce a
    lollipop $\widetilde{M}_{117}^{(5)}$, which has five edges, because deleting an edge from a graph with canonical code 95 must produce a graph with canonical code strictly less than 95, whereas
    $117>95$.
    
    \item \label{2}
    The induced subgraph $\widetilde{M}_{s}^{(b)}$ must have strictly more edges than $\widetilde{M}_{a}^{(b)}$. For example, an hourglass $\widetilde{M}_{235}^{(5)}$, which has six edges, cannot be transformed into a chevron $\widetilde{M}_{223}^{(5)}$, which has seven edges, by deleting edges, even though $s=235>a=223$.
    
    \item \label{3}
    Although a degree sequence does not in general determine a graph, among the 21 connected five-node subgraphs considered here it uniquely determines the isomorphism class except for two pairs: the banner and the lollipop, which both have degree sequence
    $(1,2,2,2,3)$, and the ufo and the house, which both have degree
    sequence $(2,2,2,3,3)$. Degree sequences can therefore provide useful intuition about which edges may be removed when moving from one subgraph to another.\\
    \\
    \noindent For a graph $H$ on $b$ nodes, let
    \[
        d_1(H)\leq d_2(H)\leq\cdots\leq d_b(H)
    \]
    denote its degree sequence in ascending order. If
    \[
        d_{\ell}(M_s^{(b)})<d_{\ell}(M_a^{(b)})
    \]
    for any $\ell=1,\ldots,b$, then the corresponding entry of $A$ is
    zero. Indeed, deleting edges cannot increase the degree of any node, and hence the ordered degree sequence of a graph obtained by edge deletion must be component-wise no greater than that of the original graph.\\
    \\
    \noindent For example, a banner, with degree sequence $(1,2,2,2,3)$, cannot be reduced to a 5-star, with degree sequence $(1,1,1,1,4)$, because its maximum degree would have to increase from 3 to 4. Similarly, an hourglass, with degree sequence $(2,2,2,2,4)$, cannot be reduced to a bull, with degree sequence $(1,1,2,3,3)$, and a chevron, with degree sequence $(2,2,2,4,4)$, cannot be reduced to a spinning top, with
    degree sequence $(1,2,3,3,3)$.
    
    \item
    Applying the necessary conditions in items~\ref{1}--\ref{3}
    identifies all zero entries of $A$ except six. These remaining cases depend on the arrangement of edges rather than only on the canonical code, number of edges, and degree sequence:
    \begin{enumerate}[label=(\alph*)]
        \item A ufo does not reduce to a bull. The ufo has degree sequence $(2,2,2,3,3)$ and the bull has degree sequence $(1,1,2,3,3)$. Since the two graphs differ by one edge, obtaining the bull degree sequence would require deleting an edge between two degree-2 nodes, but no such edge exists in the ufo.
    
        \item An hourglass reduces to a lollipop but not to a banner,
        although the lollipop and banner have the same degree sequence
        $(1,2,2,2,3)$. Every one-edge deletion from an hourglass that
        produces this degree sequence gives a lollipop rather than a
        banner.
    
        \item A ufo reduces to a banner but not to a lollipop, although these two target graphs again have the same degree sequence. In fact, deleting any one edge from a ufo gives a banner.
    
        \item A chevron reduces to a banner, through either a stingray or a ufo, but does not reduce to a lollipop.
    
        \item An arrowhead reduces to a crown or an envelope but not to a kite. Since the arrowhead and kite differ by one edge, a single edge deletion would have to create the degree-1 node of the kite. This is impossible because every node of the arrowhead has degree at least 3.
    
        \item A crown reduces to a house but not to a ufo, although the house and ufo have the same degree sequence $(2,2,2,3,3)$.
    \end{enumerate}
    
    \item For the remaining non-zero elements of $A$, the corresponding multiplicities can be calculated analytically, or evaluated directly, using \eqref{eq:m_75_5}--\eqref{eq:m_1023_5}. For example, consider the number of 5-stars $M_{75}^{(5)}$ contained in a chevron $M_{223}^{(5)}$. A chevron has degree sequence $(2,2,2,4,4)$, and \eqref{eq:m_75_5} immediately gives $\sum_{i:k_i>3}\binom{k_i}{4}=2$. Thus a chevron contains two 5-stars.\\
    \\
    \noindent A slightly more involved example is the number of crickets $M_{79}^{(5)}$ contained in the 5-complete $M_{1023}^{(5)}$. The 5-complete has degree sequence $(4,4,4,4,4)$, and \eqref{eq:m_79_5} gives
    \[
        \frac{1}{2}\sum_{i:k_i>3}(g^3)_{ii}\binom{k_i-2}{2}
        = \frac{1}{2}\tr(g^3) = 3|M_{7}^{(3)}|,
    \]
    where the last equality follows from
    \cite[eqn.~2]{agasse-duval_lawford25}. Since $K_n$ contains
    $\binom{n}{3}$ triangles, $K_5$ contains $3\binom{5}{3}=30$ crickets.\\
    \\
    \noindent The same result can also be seen directly. Each node of $K_n$ can serve as the centre of a 5-star, so the number of 5-stars in $K_n$ is $n\binom{n-1}{4} = 5\binom{n}{5}$. Hence there are five 5-stars in $K_5$. A cricket is obtained from a 5-star by adding an edge between two of its four degree-1 nodes, and
    there are $\binom{4}{2}=6$ ways to choose this edge. Thus $K_5$ contains $5\times6=30$ crickets.
    \end{enumerate}
    \end{remark}

\clearpage

\section{Alternative text for figures and graphical tables}\label{app:alt_text}
This appendix provides alternative text descriptions for the figures
and graphical tables.

\paragraph{Table~\ref{tab:subgraphs5}.}\textbf{Alt text:}
Diagrams of the 21 non-isomorphic connected simple graphs on five
nodes, arranged in three rows and seven columns. Each graph is labelled
with its descriptive name and canonical code, from the 5-star
$M_{75}^{(5)}$ to the 5-complete graph $M_{1023}^{(5)}$.

\paragraph{Table~\ref{tab:subgraphs34}.}\textbf{Alt text:}
Diagrams of the eight non-isomorphic connected simple graphs on three
or four nodes, arranged in two rows and four columns. The first row
shows the 3-star $M_{3}^{(3)}$, triangle $M_{7}^{(3)}$, 4-star
$M_{11}^{(4)}$, and 4-path $M_{13}^{(4)}$. The second row shows the
tadpole $M_{15}^{(4)}$, 4-circle $M_{30}^{(4)}$, diamond
$M_{31}^{(4)}$, and 4-complete graph $M_{63}^{(4)}$.

\paragraph{Table~\ref{tab:unwanted_cases}.}\textbf{Alt text:}
Diagram of the eight non-isomorphic connected graphs on three or four
nodes used in the five-node count formulae: the 3-star, triangle,
4-star, 4-path, tadpole, 4-circle, diamond, and 4-complete graph,
together with their canonical codes.

\paragraph{Figure~\ref{fig:six_node_graphlets}.}\textbf{Alt text:}
Two diagrams of connected six-node graphs used to illustrate
extensions of the counting strategies. Panel~(a) shows
$M_{7919}^{(6)}$, which has a unique degree-5 node whose deletion
leaves a house. Panel~(b) shows $M_{1182}^{(6)}$, used to illustrate
an extension of the walk-based counting strategy.

\paragraph{Figure~\ref{fig:illustrative_graphs}.}\textbf{Alt text:}
Three network diagrams used in the analytical examples. Panel~(a)
shows the complete graph $K_{20}$. Panel~(b) shows a complete
5-partite graph with three nodes in each part. Panel~(c) shows a
regular ring lattice with 29 nodes, in which each node is connected
to its ten nearest neighbours in each direction.

\paragraph{Figure~\ref{fig:WN_2013_4_subgraph_counts_main}.}\textbf{Alt text:} Three panels summarizing five-node subgraph counts for the Southwest Airlines network in 2013Q4. Panel~(a) compares general and induced counts for all 21 five-node subgraphs. Panel~(b) shows standardized deviations of the induced counts from an Erd{\H{o}}s--R{\'e}nyi null ensemble on a symmetric logarithmic scale; many deviations are very large and positive. Panel~(c) shows standardized deviations from a degree-preserving null ensemble; most are much closer to zero, although
several subgraphs retain deviations beyond plus or minus two.

\paragraph{Table~\ref{tab:WN_2013_4_subgraph_examples}.}\textbf{Alt text:}
Twenty-one network diagrams showing representative copies of each
connected five-node subgraph in the Southwest Airlines route network
in 2013Q4. In each panel, the five selected airports and the edges of
the corresponding subgraph are highlighted against the remainder of
the airline network, and the selected airports are labelled by their
three-letter IATA codes.

\clearpage

\singlespace

\bibliographystyle{plainnat}

\bibliography{ms}

\end{document}